\documentclass[12pt, reqno, twoside]{amsart}
\usepackage{amsmath, amsthm, mathrsfs, amscd, amsfonts, amssymb, graphicx, color}
\usepackage[bookmarksnumbered, colorlinks, plainpages]{hyperref}
\hypersetup{colorlinks=true,linkcolor=red, anchorcolor=green, citecolor=cyan, urlcolor=red, filecolor=magenta, pdftoolbar=true}
\usepackage{enumitem}
\usepackage{enumitem}
\newtheorem{theorem}{Theorem}[section]
\newtheorem{lemma}[theorem]{Lemma}

\theoremstyle{definition}

\theoremstyle{remark}
\newtheorem{remark}[theorem]{Remark}
\numberwithin{equation}{section}

\begin{document}

\setcounter{page}{1}

\title[Gradient estimates for nonlinear elliptic equations]{Souplet-Zhang and Hamilton type gradient estimates for 
nonlinear elliptic equations on smooth metric measure spaces}

\author[A. Taheri]{Ali Taheri}

\author[V. Vahidifar]{Vahideh Vahidifar}

\address{School of Mathematical and  Physical Sciences, 
University of Sussex, Falmer, Brighton, United Kingdom.}
\email{\textcolor[rgb]{0.00,0.00,0.84}{a.taheri@sussex.ac.uk}} 
 
\address{School of Mathematical and  Physical Sciences, 
University of Sussex, Falmer, Brighton, United Kingdom.}
\email{\textcolor[rgb]{0.00,0.00,0.84}{v.vahidifar@sussex.ac.uk}}



\begin{abstract}
In this article we present new gradient estimates for positive solutions to a class of nonlinear elliptic equations 
involving the $f$-Laplacian on a smooth metric measure space. The gradient estimates of interest are of 
Souplet-Zhang and Hamilton types respectively and are established under natural lower bounds on the generalised Bakry-\'Emery 
Ricci curvature tensor. From these estimates we derive amongst other things Harnack inequalities and general global 
constancy and Liouville-type theorems. The results and approach undertaken here provide a unified treatment and  
extend and improve various existing results in the literature. Some implications and applications are presented and discussed. 
\end{abstract}

\maketitle

\section{Introduction} \label{sec1}

In this paper we are concerned with deriving gradient estimates for positive solutions to a class of nonlinear 
elliptic equations on smooth metric measure spaces and exploiting some of their local and global implications. 
In more specific terms, we aim to formulate and prove gradient estimates of Souplet-Zhang and Hamilton types 
for positive smooth solutions $u$ to the nonlinear elliptic equation: 
\begin{align} \label{eq11}
\Delta_f u + \Sigma (x, u) =0,  \qquad \Sigma: M \times \mathbb{R} \to \mathbb{R},  
\end{align}
on a smooth metric measure space $(M, g, d\mu)$, where $M$ is a Riemannian manifold with dimension 
$n \ge 2$, $d\mu=e^{-f} dv_g$ is a weighted measure, $f$ is a smooth potential on $M$, $g$ is the 
Riemannian metric and $dv_g$ is the usual Riemannain volume measure. 
Here $\Delta_f$ is the so-called $f$-Laplacian that acts on functions $u \in \mathscr{C}^2(M)$ by 
$\Delta_f u = \Delta u - \langle \nabla f, \nabla u\rangle$, and $\Sigma=\Sigma(x,u)$ 
is a sufficiently smooth but otherwise arbitrary {\it nonlinearity} depending on the spatial variable $x$ and the 
independent variable $u$. 

The $f$-Laplacian appearing in \eqref{eq11} is a symmetric diffusion operator with respect to the invariant weighted 
measure $d \mu =e^{-f} dv_g$. It arises naturally in a variety of contexts ranging from geometry, probability theory and 
stochastic processes to quantum field theory, statistical mechanics and kinetic theory \cite{Bak, LiX, VC, Wang, WuLM}. 
It is a generalisation of the Laplace-Beltrami operator to the smooth metric measure space context and 
one recovers the latter exactly when the potential $f$ is a constant. 

Alternatively and more directly one can view the $f$-Laplacian as the infinitesimal generator of the Dirichlet form 
\begin{equation}
(u,v) \mapsto \mathscr{E}[u,v] = \int_M \langle \nabla u, \nabla v \rangle \, e^{-f} dv_g(x), 
\qquad \forall u,v \in \mathscr{C}_c^\infty(M). 
\end{equation}

Here Ito's SDE theory or Dirichlet form theory imply the existence of a minimal diffusion process 
$(X_t:t<T)$ on $M$ with infinitesimal generator $L=\Delta_f=e^f {\rm div} (e^{-f} \nabla)$ and 
lifetime $T$. This diffusion process is the solution to the Ito-Langevin SDE 
\begin{equation}
dX_t = \sqrt{2} dW_t - \nabla f (X_t) \, dt, 
\end{equation}
where $W_t$ denotes the standard Riemannian Brownian motion on $M$ and so as a result the $f$-heat equation 
$\partial_t u = \Delta u - \langle \nabla f, \nabla u \rangle$ is then the backward Kolmogorov equation associated with 
$(X_t: t<T)$. Note in particular that the $f$-heat semigroup $e^{t \Delta_f}$ here can be represented in terms 
of the expectation 
\begin{equation}
e^{t \Delta_f} v (x) = \mathbb{E}_x [v(X_t) 1_{\{t<T\}}], \qquad \forall v \in \mathscr{C}_c^\infty(M),
\end{equation}
and the equilibrium solutions to the $f$-heat equation $\partial_t u = \Delta u - \langle \nabla f, \nabla u \rangle$ 
are precisely the $f$-harmonic functions $\Delta_f u =0$, the linear version of \eqref{eq11}. 

Apart from the connection alluded to above there is also a close connection between symmetric 
diffusion operators and Schr\"odinger operators 
\cite{LiX, WuLM}. Indeed the operator $L=\Delta_f$ viewed in $L^2(M, d\mu)$ is unitarily 
equivalent to the Schr\"odinger operator $H=\Delta - V$ with $V=|\nabla f|^2/4 - \Delta f/2$ viewed in $L^2(M, dv_g)$ through the 
unitary isomorphism $\mathscr L$ given by 
\begin{align}
\mathscr{L}: &L^2(M, d\mu) \to L^2(M, dv_g), \nonumber \\
&u \mapsto \mathscr{L} u = e^{-f/2} u.
\end{align} 
It is easily seen that 
\begin{equation}
H(\mathscr{L} u)=[\Delta - |\nabla f|^2/4 + \Delta f/2](e^{-f/2} u) = e^{-f/2} \Delta_f u = \mathscr{L} \Delta_f u,
\end{equation} 
and so in particular $\Delta_f u = 0 \iff H (e^{-f/2} u) =0$; giving a correspondence between $f$-harmonic 
functions and the solutions to the stationary Schr\"odinger equation. Moreover, on the level of the heat 
semigroups we have the relation 
\begin{equation}
e^{tH} = \mathscr{L} \circ e^{t \Delta_f} \circ \mathscr{L}^{-1},  \qquad t \ge 0, 
\end{equation}
where by the Feynman-Kac formula and upon recalling $H=\Delta - V$ with $V=|\nabla f|^2/4 - \Delta f/2$ 
from above we have for all $v \in \mathscr{C}_c^\infty(M)$: 
\begin{equation}
e^{tH} v (x) 
= \mathbb{E}_x \left[ v(W_t) {\rm exp} \left\{ - \int_0^t \left( \frac{|\nabla f(W_t)|^2}{4} 
- \frac{\Delta f (W_s)}{2} \right) \, ds \right\} 1_{\{t<T\}} \right].  
\end{equation}

This results in a correspondence between the $L^p$-uniqueness of $f$-heat semigroup 
$e^{t \Delta_f}$ and the Schr\"odinger semigroup $e^{tH}$ which is fundamental in the 
well-posedness of their respective Cauchy problems. 
There are other important and close connections between these semigroups and 
evolution equations with one another and, e.g., the Fokker-Planck equation (the forward Kolmogorov 
equation associated with the diffusion process $(X_t : t<T)$ above) that we do not 
discuss here. (See \cite{Bak, LiX, Gross, VC, Wang, WuLM} and 
the references therein for further detail and discussion in this direction.)

Returning to the smooth metric measure space context and attempting to describe the geometry 
of the $f$-Laplacian and its associated diffusion process, we next introduce the generalised 
Ricci curvature tensor for the triple $(M, g, d\mu)$ by setting
\begin{equation} \label{Ricci-m-f-intro}
{\mathscr Ric}_f(g) = {\mathscr Ric}(g) + {\rm Hess}(f).   
\end{equation}
Here ${\mathscr Ric}(g)$ denotes the usual Riemannain Ricci curvature tensor of $g$ [{\it see} \eqref{RicciIndex-def}] and 
${\rm Hess}(f)=\nabla \nabla f$ stands for the Hessian of $f$. 
The weighted Bochner-Weitzenb\"ock formula relating the quantities $\Delta_f |\nabla u|^2$, $\Delta_f u$ 
and ${\mathscr Ric}_f(\nabla u, \nabla u)$ for any $u \in {\mathscr C}^3(M)$ then asserts that 
\begin{equation} \label{Bochner}
\frac{1}{2} \Delta_f |\nabla u|^2 = |{\rm Hess}(u)|^2 + \langle \nabla u, \nabla \Delta_f u \rangle 
+ {\mathscr Ric}_f (\nabla u, \nabla u),  
\end{equation}
with the norm of the Hessian on the right being the Hilbert-Schmidt norm of $2$-tensors. 
As a result of \eqref{Bochner} a curvature lower bound ${\mathscr Ric}_f(g) \ge \mathsf{k} g$ together with the Cauchy-Schwarz inequality 
result in the curvature-dimension condition ${\rm CD}(\mathsf{k}, \infty)$, that is, 
\footnote{As the right-hand side of \eqref{eq-CD-condition} contains $\Delta u$ and not $\Delta_f u$, 
in general, one can only have ${\rm CD}(\mathsf{k}, \infty)$. If $f$ is constant, then $\Delta_f=\Delta$, 
${\mathscr Ric}_f(g) = {\mathscr Ric}(g)$ and \eqref{eq-CD-condition} gives ${\rm CD}(\mathsf{k}, n)$. 
See \cite{Bak} for more discussion on the significance of the curvature-dimension condition 
${\rm CD}({\mathsf k}, m)$ in analysis and geometry.
}
\begin{align} \label{eq-CD-condition}
\frac{1}{2} \Delta_f |\nabla u|^2 - \langle \nabla u, \nabla \Delta_f u \rangle 
&= |{\rm Hess}(u)|^2 + {\mathscr Ric}_f (\nabla u, \nabla u) 
\ge \frac{(\Delta u)^2}{n} + \mathsf{k} |\nabla u|^2.
\end{align}

The main objective in this paper is to develop local and global gradient estimates of Souplet-Zhang 
and Hamilton types along with Harnack inequalities for positive smooth solutions to \eqref{eq11}. It 
is well-known that such estimates and inequalities form the basis for deriving various qualitative 
properties of solutions and are thus of great significance \cite{Gr, LY86, LiP-book, SchYau, YauH, Zhang}. 
Such properties may include local and global bounds, eigenvalues and spectral asymptotics, Liouville-type 
results, heat kernel bounds and asymptotics, applications to the time evolution equations and/or geometric 
flows, characterisation of ancient and eternal solutions, analysis of singularities and much more ({\it see} 
\cite{AM, Bid, Bri, CGS, ChYau, Dung, Giaq, Gr, Ham, KT, LiP-book, LiX, PG, SZ, Taheri-book-one, Taheri-book-two, 
Taheri-GE-1, Taheri-GE-2, TVahGrad, TahVah, Wang, Wu15, Wu18} and the references therein for further details and 
discussion).

Whilst in formulating and proving gradient estimates one often works with a specific choice or type of nonlinearity with 
an explicit structure of singularity, regularity, growth and decay, in this article we keep the analysis and discussion on 
a fairly general level without confining to specific examples or choices. This way we firstly provide a unified treatment 
of the estimates and secondly examine more transparently how the structure and form of the nonlinearity ultimately 
influences and interacts with the estimates and the subsequent results. As such the approach and analysis 
undertaken here largely unify, extend and improve various existing results in the literature for specific choices 
of nonlinearities. We discuss this further in the overview below.

Gradient estimates for positive solutions to \eqref{eq11} with logarithmic type nonlinearities in the form:   
\begin{equation} \label{eq1.4}
\Delta_f u + \mathsf{A}(x) u \log u + \mathsf{B}(x) u = 0, 
\end{equation}
along with their parabolic counterparts have been the subject of extensive studies 
({\it see}, e.g., \cite{Ma, Ruan, YYY} and the references therein). The interest in such 
problems originate partly from its natural links with gradient Ricci solitons and partly from 
geometric and functional inequalities, most notably, the logarithmic Sobolev inequalities \cite{Bak, Gross, VC, Zhang}. 
Recall that a Riemannian manifold 
$(M,g)$ is said to be a gradient Ricci soliton {\it iff} there there exists a smooth function $f$ on $M$ and 
a constant $\mu \in \mathbb{R}$ such that ({\it cf.} \cite{Cao, Chow, Lott})
\begin{equation}\label{eq1.5}
{\mathscr Ric}_f(g) = {\mathscr Ric}(g) + {\rm Hess}(f) = \mu g.
\end{equation}
The notion is clearly an extension of an Einstein manifold and has a fundamental role in the analysis of singularities 
of the Ricci flow \cite{Ham, Zhang}. Taking trace from both sides of \eqref{eq1.5} and using the contracted 
Bianchi identity leads one to a simple form of \eqref{eq1.4} with constant coefficients: 
$\Delta u + 2 \mu u \log u = \lambda u$ for suitable constant $\lambda$ and $u=e^f$ 
({\it see} \cite{Ma} for details). Other types of equations generalising and strengthening \eqref{eq1.4} are 
\begin{equation}
\Delta_f u + \mathsf{A}(x) u^s \gamma(\log u) + \mathsf{B}(x) u^p + \mathsf{C}(x) u^q=0, 
\qquad \gamma \in \mathscr{C}^1(\mathbb{R}), 
\end{equation}
with real exponents $r, p, q$ and have been studied in detail in \cite{CGS, GKS, Taheri-GE-1, Taheri-GE-2, Wu15}.

Yamabe type equations $\Delta u + \mathsf{A}(x) u^p + \mathsf{B}(x) u =0$ also come under the form  
\eqref{eq11} with a power-like nonlinearity. In \cite{Bid} 
the equation $\Delta u + u^p+ \mathsf{B} u =0$ is studied on a compact manifold and under suitable 
conditions on ${\mathscr Ric}(g)$, $n$, $p$ and $\mathsf{B}$ it is shown to admit only constant 
solutions. The equation $\Delta u + \mathsf{A}(x) u^p =0$ for $1 \le p < 2^\star-1=(n+2)/(n-2)$ with $n \ge 3$
is considered in \cite{GidSp} and it is proved that when ${\mathscr Ric}(g) \ge 0$ any non-negative 
solution to this equation must be zero. In \cite{Yang} it is shown that the same equation with constant 
$\mathsf{A}>0$, $p<0$ admits no positive solution when ${\mathscr Ric}(g) \ge 0$. For a detailed 
account on the Yamabe problem in geometry see \cite{Lee, Mast}. The natural form of Yamabe 
equation in the setting of smooth metric measure spaces is   
\begin{equation} 
\Delta_f u + \mathsf{A}(x) u^p+ \mathsf{B}(x) u  = 0. 
\end{equation}
For gradient estimates, Harnack inequalities and other counterparts of the above results we 
refer the reader to \cite{Case, Wu15, ZhMa}.

A far more general form of Yamabe equation is the Einstein-scalar field Lichnerowicz equation 
(see, e.g., \cite{CBY, Chow}). 
Here when the underlying manifold has dimension $n \ge 3$ this takes the form 
$\Delta u + \mathsf{A}(x)u^p + \mathsf{B}(x) u^q + \mathsf{C}(x) u =0$ with 
$p=(n+2)/(n-2)$ and $q=(3n-2)/(n-2)$ while when $n=2$ this takes the 
form $\Delta u + \mathsf{A}(x) e^{2u} + \mathsf{B}(x) e^{-2u} + \mathsf{C}(x)=0$.  
The Einstein-scalar field Lichnerowicz equation in the setting of smooth metric 
measure spaces can naturally be formulated as:  
\begin{equation}
\Delta_ f u + \mathsf{A} u^p + \mathsf{B} u^q + \mathsf{C} u \log u + \mathsf{D} u =0, 
\end{equation}
and 
\begin{equation}
\Delta_ f u + \mathsf{A} e^{2u} + \mathsf{B} e^{-2u} + \mathsf{C} =0,  
\end{equation}
with $\mathsf{A}=\mathsf{A}(x)$, $\mathsf{B}=\mathsf{B}(x)$, $\mathsf{C}=\mathsf{C}(x)$ and $\mathsf{D}=\mathsf{D}(x)$.

For gradient estimates, Harnack inequalities and Liouville-type results in this and related contexts we refer the reader to 
\cite{Dung, Taheri-GE-1, Taheri-GE-2, Wu18} and the references therein.

\qquad \\
{\bf Plan of the paper.} Let us bring this introduction to an end by briefly describing the plan of the paper. 
In the concluding part of Section \ref{sec1} (below) we describe the notation and terminology as used in the paper. 
In Section \ref{sec2} we present the main results of the paper, specifically, a local and global gradient estimate 
of Souplet-Zhang type for equation \eqref{eq11}, followed by local and global Harnack inequalities, 
a local and global gradient estimate of Hamilton type for equation \eqref{eq11} and a general Liouville-type 
result with some concrete applications to particular classes of nonlinearities. The subsequent sections are 
then devoted to the detailed proofs respectively. Section \ref{sec3} is devoted to the proof of Theorem 
\ref{thm1} and Section \ref{sec4} presents the proof of the local and global Harnack inequalities in 
Theorem \ref{cor Harnack}. Section \ref{sec5} is devoted to the proof of Theorem \ref{thm18} and the 
Liouville-type result in Theorem \ref{thm Liouville}, which is an application of the former, is proved in 
Section \ref{Liouville-Section}. Finally in Section \ref{Hamilton-Section} we present the proof of the 
global Hamilton bound in Theorem \ref{global-Hamilton-Thm}. We also discuss some applications 
and examples.

\qquad \\
{\bf Notation.} 
Fixing a reference point $p \in M$ we denote by $d=d_p(x)$ the Riemannian distance between $x$ 
and $p$ and by $r=r_p(x)$ the geodesic radial variable with origin at $p$. We 
write $\mathcal{B}_a(p) \subset M$ for the geodesic ball of radius $a>0$ centred at $p$. 
When the choice of the point $p$ is clear from the context we often abbreviate and write 
$d(x)$, $r(x)$ or $\mathcal{B}_a$ respectively. 
For $X \in {\mathbb R}$ we write $X_+=\max(X, 0)$ and $X_-=\min(X, 0)$. 
Hence $X=X_+ + X_-$ with $X_+ \ge 0$ and $X_- \le 0$. 
For the sake of future reference and use in the gradient estimates to come, we also introduce the 
functional quantity    
\begin{equation} \label{sigma def}
[\gamma_{\Delta_f}]_+=\max(\gamma_{\Delta_f}, 0), 
\end{equation}
where we have set
\begin{equation} \label{sigma def alt}
\gamma_{\Delta_f} = \max_{\partial {\mathcal B}_1(p)} \Delta_f r  = \max \{ \Delta_f r (x) : d_p(x) =1 \}, 
\end{equation}
that is, the maximum of the $f$-Laplacian of the radial variable $r=r_p(x)$, that is, $\Delta_f r$, taken over the unit 
sphere centred at the reference point $p$.

We typically denote partial derivatives with subscripts, in particular, for the nonlinear function $\Sigma=\Sigma(x,u)$ 
we make use of $\Sigma_x$, $\Sigma_u$, {\it etc.} and we reserve the notation $\Sigma^x$ for 
the function $\Sigma(\cdot, u)$ obtained by freezing the argument $u$ and viewing it as 
a function of $x$ only. In fact below we frequently speak of $\nabla \Sigma^x$, 
$\Delta \Sigma^x$ and $\Delta_f \Sigma^x$.   

For the sake of reader's convenience we recall that in local coordinates $(x^i)$ we have the following 
formulae for the Laplace-Beltrami operator, Riemann and Ricci curvature tensors respectively:  
\begin{equation} \label{Lap-coordinate}
\Delta = \frac{1}{\sqrt{|g|}} \frac{\partial}{\partial x_i} \left( \sqrt{|g|} g^{ij} \frac{\partial}{\partial x_j} \right),
\end{equation}
and 
 \begin{equation}
[{\mathscr Rm}(g)]^\ell_{ijk} = \frac{\partial \Gamma^\ell_{jk}}{\partial x_i} - \frac{\partial \Gamma^\ell_{ik}}{\partial x_j} 
+ \Gamma^p_{jk} \Gamma^\ell_{ip}  - \Gamma^p_{ik} \Gamma^\ell_{jp},   
\end{equation}
and 
\begin{equation} \label{RicciIndex-def}
[{\mathscr Ric}(g)]_{ij} = \frac{\partial \Gamma_{ij}^k}{\partial x_k} - \frac{\partial \Gamma_{\ell j}^\ell}{\partial x_i} 
+ \Gamma^k_{ij} \Gamma^\ell_{\ell k} - \Gamma^\ell_{ik} \Gamma^k_{\ell j}. 
\end{equation}
We point out that here
\begin{equation}
\Gamma^k_{ij} = \frac{1}{2} g^{k\ell} \left( \frac{\partial g_{j \ell}}{\partial x_i} 
+ \frac{\partial g_{i \ell}}{\partial x_j} - \frac{\partial g_{ij}}{\partial x_\ell} \right),
\end{equation}
are the Christoffel symbols. Additionally $g_{ij}$, $|g|$ and $g^{ij} = (g^{-1})_{ij}$ are the components, 
the determinant and the components of the inverse of the metric tensor $g$ respectively. Note in 
particular that upon referring to \eqref{Lap-coordinate} we also have 
\begin{equation}
\Delta_f = \frac{1}{\sqrt{|g|}} \frac{\partial}{\partial x_i} \left( \sqrt{|g|} g^{ij} \frac{\partial}{\partial x_j} \right) 
- g^{ij} \frac{\partial f}{\partial x_i} \frac{\partial}{\partial x_j},
\end{equation}
for the $f$-Laplacian.

\section{Statement of the main results} \label{sec2}

In this section we present the main results of the paper with some accompanying discussion. The detailed arguments 
and proofs are delegated to the subsequent sections. We begin with a local Souplet-Zhang type gradient estimates for 
\eqref{eq11} exploiting the role of the nonlinearity $\Sigma$ in \eqref{eq11}.

An important quantity appearing 
in this estimate [{\it cf.} \eqref{eq13} below] is $\mathsf{R}^{1/2}_\Sigma(u) + \mathsf{P}^{1/3}_\Sigma(u)$ 
[see \eqref{eq2.2}-\eqref{eq2.3} for the explicit formulations].  As a close inspection reveals, these terms 
directly link to the nonlinearity $\Sigma$ and the solution $u$, and as will be seen later, they play a decisive 
role in global estimates, Harnack inequalities and the Liouville theorems that follow. Note that since $u>0$ 
and $\mathcal{B}_{R} \subset M$ is compact, $u$ is bounded away from zero and bounded from above; 
hence, in particular, this quantity is finite. The other quantities appearing here are $k$ and 
$[\gamma_{\Delta_f}]_+$ relating to the curvature bound ${\mathscr Ric}_f(g) \ge -(n-1) k g$ 
in the theorem, the potential $f$, and the geometry of the triple $(M, g, d\mu)$ [{\it cf.} \eqref{sigma def}].

\begin{theorem} \label{thm1}
Let $(M, g, e^{-f} dv_g)$ be a complete smooth metric measure space satisfying ${\mathscr Ric}_f (g) \ge - (n-1) k g$ 
in $\mathcal{B}_{R}=\mathcal{B}_R(p) \subset M$ where $R \ge 2$ and $k \ge 0$. Let $u$ be a positive solution to 
the nonlinear elliptic equation \eqref{eq11} in $\mathcal{B}_R$ with $0<u \le A$. Then on $\mathcal{B}_{R/2}$ 
the solution $u$ satisfies the estimate:
\begin{align} \label{eq13}
\frac{|\nabla u|}{u} \le   
C \left\{ \frac{1}{R} +  \sqrt{k}  + \sqrt{ \frac{[\gamma_{\Delta_f}]_+}{R}}
+ \sup_{\mathcal{B}_{R}} \left[ \mathsf{R}^{1/2}_\Sigma(u) 
+ \mathsf{P}^{1/3}_\Sigma(u) \right] \right\} \left(1 - \log \frac{u}{A} \right). 
\end{align}
Here $C>0$ is a constant depending only on $n$ and the quantities $\mathsf{P}_\Sigma(u)$ and $\mathsf{R}_\Sigma(u)$ 
are defined respectively by 
\begin{equation} \label{eq2.2}
\mathsf{P}_\Sigma(u) = \frac{|\Sigma_x(x,u)|}{u [1- \log(u/A)]^2}, 
\end{equation} 
and
\begin{equation} \label{eq2.3}
\mathsf{R}_\Sigma(u) = \left[ \frac{u[1-\log(u/A)] \Sigma_u(x,u) + \log(u/A) \Sigma (x,u)}{u [1-\log(u/A)]^2} \right]_+,  
\end{equation}
whilst the functional quantity $[\gamma_{\Delta_f}]_+ \ge 0$ is as in \eqref{sigma def}. 
\end{theorem}

The local estimate above has a global counterpart subject to the prescribed bounds in the theorem being global. 
The proof follows by passing to the limit $R \to \infty$ in \eqref{eq13}. The precise 
formulation of this is given in the following theorem.

\begin{theorem} \label{thm1-global}
Let $(M, g, e^{-f} dv_g)$ be a complete smooth metric measure space satisfying the global curvature lower bound 
${\mathscr Ric}_f (g) \ge -(n-1) k g$  on $M$ with $k \ge 0$. Suppose $u$ is a positive solution to the nonlinear 
elliptic equation \eqref{eq11} in $M$ satisfying $0<u \le A$. Then $u$ satisfies the following global estimate on $M$:
\begin{align} \label{eq13-global}
\frac{|\nabla u|}{u} \le C \left\{ \sqrt{k} 
+ \sup_{M} \left[ \mathsf{R}^{1/2}_\Sigma(u) 
+ \mathsf{P}^{1/3}_\Sigma(u) \right] \right\} \left(1 - \log \frac{u}{A} \right). 
\end{align}
Here $C>0$ depends only on $n$, 
the quantities $\mathsf{R}_{\Sigma}(u)$ and $\mathsf{P}_{\Sigma}(u)$ are as in \eqref{eq2.2} and \eqref{eq2.3} in 
Theorem $\ref{thm1}$ and the supremum in \eqref{eq13-global} is taken over $M$.
\end{theorem}

One of the useful consequences of the estimates above is the following elliptic Harnack inequality for bounded 
positive solutions to the nonlinear elliptic equation \eqref{eq11}. Later on we will prove another version of this inequality 
using a different approach.

\begin{theorem} \label{cor Harnack}
Under the assumptions of Theorem $\ref{thm1}$, if $u$ is a positive solution to the nonlinear elliptic equation \eqref{eq11} 
with $0<u \le A$, then for all $x$, $y$ in $\mathcal{B}_{R/2}$ we have:
\begin{equation} \label{eq2.5}
0< u(x) \le (eA)^{1-\gamma} u^\gamma(y), 
\end{equation}
where the exponent $\gamma \in (0,1)$ in \eqref{eq2.5} is given explicitly by the exponential 
\begin{equation} \label{exponent-local}
\gamma = {\rm exp} \left[ - C d(x, y) 
\left( \frac{1}{R} +  \sqrt{k}  + \sqrt{ \frac{[\gamma_{\Delta_f}]_+}{R}} 
+ \sup_{\mathcal{B}_R} \left[ \mathsf{R}^{1/2}_\Sigma(u) + \mathsf{P}^{1/3}_\Sigma(u) \right] \right) \right].
\end{equation}
Additionally, subject to the prescribed global bounds in Theorem $\ref{thm1-global}$, for all $x$, $y$ 
in $M$, we have the same inequality \eqref{eq2.5} now with the exponent
\begin{equation} \label{exponent-global}
\gamma = {\rm exp} \left[ - C d(x, y) 
\left( \sqrt k + \sup_{M} \left[ \mathsf{R}^{1/2}_\Sigma(u) + \mathsf{P}^{1/3}_\Sigma(u) \right] \right) \right].
\end{equation}
In \eqref{exponent-local} and \eqref{exponent-global}, $d(x,y)$ is the geodesic distance between $x$ and $y$ 
and $C>0$ is as in \eqref{eq13} and \eqref{eq13-global} respectively. 
\end{theorem}

The next type of estimate we move on to is a local Hamilton-type gradient estimate for positive solutions 
to the nonlinear elliptic equation \eqref{eq11}. Here the estimate makes use of two non-negative 
parameters $\alpha$, $\beta$ with appropriate ranges ({\it see} below for details). This gives the estimate 
flexibility and scope for later applications.

The important quantity appearing in the estimate \eqref{eq1.12} here is 
$\mathsf{T}_\Sigma^{1/2} (u) + \mathsf{S}_\Sigma^{1/3} (u)$ [see \eqref{eq2.10}-\eqref{eq2.10-alt} 
for the explicit description of the terms] which again directly links to the nonlinearity $\Sigma$, the 
parameters $\alpha$, $\beta$ and the solution $u$. The other quantities appearing here are $k$ 
and $[\gamma_{\Delta_f}]_+$ relating to the curvature bound ${\mathscr Ric}_f(g) \ge -(n-1) k g$ 
in the theorem, the potential $f$, and the geometry of the triple $(M, g, d\mu)$.

\begin{theorem} \label{thm18}
Let $(M, g, e^{-f} dv_g)$ be a complete smooth metric measure space satisfying ${\mathscr Ric}_f (g) \ge -(n-1) k g$ 
in $\mathcal{B}_{R}$ with $k\ge 0$. Let $u$ be a positive solution to the nonlinear elliptic equation \eqref{eq11}. 
Then for every $\alpha \ge 0$, $0<\beta<1/(1+\alpha)$ the solution $u$ satisfies on $\mathcal{B}_{R/2}$ the estimate:
\begin{align} \label{eq1.12}
\frac{|\nabla u|}{u^{1-\beta(\alpha/2+1)}} \le C \left\{ \frac{1}{R} + \sqrt k + \sqrt{\frac{[\gamma_{\Delta_f}]_+}{R}} 
+ \sup_{\mathcal{B}_{R}} \left[ \mathsf{T}_\Sigma^{1/2} (u)+ \mathsf{S}_\Sigma^{1/3} (u) \right] \right\} 
\Big( \sup_{\mathcal{B}_{R}} u^{\beta(\alpha/2+1)} \Big).  
\end{align}
Here $C>0$ is a constant depending only on $n$, $\alpha$, $\beta$ and the quantities 
$\mathsf{S}_{\Sigma}(u)$ and $\mathsf{T}_{\Sigma}(u)$ are defined respectively by 
\begin{equation} \label{eq2.10}
\mathsf{S}_\Sigma(u) = \frac{|\Sigma_x (x,u)|}{u}, 
\end{equation}
and
\begin{equation} \label{eq2.10-alt}
\mathsf{T}_\Sigma (u) = 2 \left[ \frac{u \Sigma_u(x,u)-[1-\beta(\alpha/2+1)]\Sigma(x,u)}{u} \right]_+,   
\end{equation}
whilst the functional quantity $[\gamma_{\Delta_f}]_+ \ge 0$ is as in \eqref{sigma def}. 
\end{theorem}

Again, the local estimate above has a global counterpart, when the asserted bounds in the theorem are global. 
This is the content of the following theorem.

\begin{theorem} \label{thm18-global}
Let $(M, g, e^{-f} dv_g)$ be a complete smooth metric measure space satisfying the global curvature lower bound 
${\mathscr Ric}_f (g) \ge -(n-1) k g$  on $M$ with $k \ge 0$. Assume $u$ is a positive solution to the nonlinear 
elliptic equation \eqref{eq11}. Then for every $\alpha \ge 0$ and $0<\beta<1/(1+\alpha)$, $u$ satisfies the following 
global estimate on $M$:
\begin{align} \label{eq1.12-global}
\frac{|\nabla u|}{u^{1-\beta(\alpha/2+1)}} 
\le C \left\{ \sqrt k + \sup_{M} \left[ \mathsf{T}_\Sigma^{1/2} (u) + \mathsf{S}_\Sigma^{1/3} (u) \right] \right\}
\Big( \sup_{M} u^{\beta(\alpha/2+1)} \Big). 
\end{align}
Here $C>0$ depends only on $n$, $\alpha$, $\beta$,   
the quantities $\mathsf{S}_{\Sigma}(u)$ and $\mathsf{T}_{\Sigma}(u)$ 
are as in \eqref{eq2.10} and \eqref{eq2.10-alt} in Theorem $\ref{thm18}$ 
and the supremums in \eqref{eq1.12-global} are taken over $M$.
\end{theorem}

\begin{remark} \label{other-curvature-assuptions-remark} {
Before moving on to discussing some applications of the above results, let us pause briefly to make a comment 
on the form of the local estimates \eqref{eq13} and \eqref{eq1.12} under a different but closely related 
curvature condition. Towards this end, recall that the generalised Bakry-\'Emery $m$-Ricci curvature 
tensor associated with the $f$-Laplacian is defined for the range $m \ge n$ by the symmetric $2$-tensor 
\begin{equation} \label{Ric f m}
{\mathscr Ric}_f^m (g) = {\mathscr Ric}(g) + {\rm Hess}(f) - \frac{\nabla f \otimes \nabla f}{m-n}. 
\end{equation}
When $m=n$, to make sense of \eqref{Ric f m} one only allows constant functions $f$ as admissible 
potentials, in which case $\Delta_f=\Delta$ and ${\mathscr Ric}_f^n (g) \equiv {\mathscr Ric} (g)$, 
and when $m=\infty$, by formally taking the limit $m \to \infty$ in \eqref{Ric f m} one sets 
\begin{equation} \label{Ric f infty}
{\mathscr Ric}^\infty_f (g) := {\mathscr Ric}_f(g) = {\mathscr Ric}(g) + {\rm Hess}(f).
\end{equation}

Referring to \eqref{Ric f m} and \eqref{Ric f infty} now, it is easily seen that a lower 
bound on ${\mathscr Ric}_f^m(g)$ is a stronger condition than a lower bound 
on ${\mathscr Ric}_f(g)$ as the former implies the latter but not {\it vice versa}. 
It can be shown [{\it see} Remark \ref{Ricmf-proof-1-remark} and 
Remark \ref{Ricmf-proof-2-remark}] that if in either of Theorem \ref{thm1} 
or Theorem \ref{thm18} the condition ${\mathscr Ric}_f(g) \ge -(n-1) kg$ 
is replaced by the {\it stronger} one ${\mathscr Ric}_f^m(g) \ge -(m-1) kg$ 
for some $n \le m <\infty$ and $k \ge 0$ then in the local estimates 
\eqref{eq13} or \eqref{eq1.12} the term containing the functional quantity 
\eqref{sigma def}, namely, the quotient 
\begin{equation} \label{extra-term-Ricmf}
[\gamma_{\Delta_f}]_+^{1/2}/\sqrt R,
\end{equation} 
can be removed from the right-hand sides respectively. This means that (upon adjusting 
the constant $C>0$ which will now depend on $m$ as well) the local estimates respectively 
take the the forms 
\begin{align} \label{eq13-Ricfm}
\frac{|\nabla u|}{u} \le   
C \left\{ \frac{1}{R} +  \sqrt{k} + \sup_{\mathcal{B}_{R}} \left[ \mathsf{R}^{1/2}_\Sigma(u) 
+ \mathsf{P}^{1/3}_\Sigma(u) \right] \right\} \left(1 - \log \frac{u}{A} \right),  
\end{align}
and
\begin{align} \label{eq1.12-Ricfm}
\frac{|\nabla u|}{u^{1-\beta(\alpha/2+1)}} \le C \left\{ \frac{1}{R} 
+ \sqrt k + \sup_{\mathcal{B}_{R}} \left[ \mathsf{T}_\Sigma^{1/2} (u) 
+ \mathsf{S}_\Sigma^{1/3} (u) \right] \right\} 
\Big( \sup_{\mathcal{B}_{R}} u^{\beta(\alpha/2+1)} \Big).   
\end{align}

As is readily seen this stronger curvature condition leads to a potentially faster decay for the 
$R$-dependent terms on the right-hand side of the estimates as $R \to \infty$. Note also that 
for $m=n$ where one recovers the Laplace-Beltrami operator with the Riemannian 
volume measure and the usual Ricci curvature lower bound ${\mathscr Ric}(g) \ge -(n-1) kg$, 
our results and estimates for the nonlinear equation \eqref{eq11} are, to 
the best of our knowledge, new in this generality for $\Sigma=\Sigma(x,u)$. 
}
\end{remark}

Let us now move on to discussing some applications of the above estimates. 
Here our focus will be on Liouville theorems and we start with a brief outline 
for such results in the linear case, namely, for $f$-harmonic 
functions under various curvature conditions, before moving 
to the full nonlinear case in \eqref{eq11}.

To this end recall that a classical form of Liouville's theorem for harmonic functions asserts that 
any positive harmonic function on $M={\mathbb R}^n$ must be a constant. (Hence by linearity 
and translation also any harmonic function on ${\mathbb R}^n$ that is bounded either from 
above or below.) The assertion is a straightforward consequence of the mean-value property: 
pick $x \neq y$ in ${\mathbb R}^n$, $R>0$ sufficiently large and $r=R-|x-y|$. Then, by the 
positivity of $u$ and the inclusion ${\mathcal B}_r(y) \subset {\mathcal B}_R(x)$ it is plain that  
\begin{equation*}
u(y) = \frac{1}{|{\mathcal B}_r(y)|} \int_{{\mathcal B}_r(y)} u 
\le \frac{|{\mathcal B}_R(x)|}{|{\mathcal B}_r(y)|} 
\left( \frac{1}{|{\mathcal B}_R(x)|} \int_{{\mathcal B}_R(x)} u \right) 
= \frac{|{\mathcal B}_R(x)|}{|{\mathcal B}_r(y)|} u(x). 
\end{equation*}
Passing to the limit $R \nearrow \infty$ (and noting $R/r \searrow 1$) gives $u(y) \le u(x)$. Interchanging 
the roles of $x$ and $y$ gives $u(x) \le u(y)$. Thus putting these together gives $u(x)=u(y)$. The arbitrariness of $x$, $y$ now 
proves the assertion.

Using gradient estimates the same Liouville property was proved by S.T.~Yau in \cite{YauH} for 
positive harmonic functions on a complete Riemennian manifold with non-negative Ricci curvature 
(see also Theorem 3.1 and Corollary 3.1 in \cite{SchYau}). An essentially similar type of gradient 
estimate can be used to prove that any positive $f$-harmonic function $u$ 
(i.e., $\Delta_f u=0$) on a complete smooth metric measure space with ${\mathscr Ric}_f^m(g) \ge 0$ for some 
$n \le m <\infty$ must be a constant ({\it see}, e.g., \cite{LiX} Theorem 1.3). Interestingly, however, the same 
conclusion dramatically fails under ${\mathscr Ric}_f(g) \ge 0$ 
alone as easy examples suggest. (Take $M={\mathbb R}^n$, $f= \langle \alpha, x \rangle$ and 
$u={\rm exp}(\langle \alpha, x \rangle)$ with $\alpha$ a fixed non-zero vector. 
Then $u>0$, $\Delta_f u = |\alpha|^2 u - \langle \alpha, \nabla u \rangle =0$ and 
${\mathscr Ric}_f(g)={\rm Hess}(f)=0$ yet $u$ is non-constant!) 
Notice that in this example ${\mathscr Ric}_f^m(g) = - [\alpha \otimes \alpha]/(m-n) \not\ge 0$ 
[{\it cf}. \eqref{Ric f m}] and that $u$ is unbounded. It is proved in \cite{Bri} that if $u$ is a {\it bounded} $f$-harmonic function 
on $M$ (e.g., a positive and bounded from above $f$-harmonic function) and ${\mathscr Ric}_f(g) \ge 0$ 
then $u$ must be a constant.

We can now formulate the following theorem which is the counterpart of the above 
for the full nonlinear equation \eqref{eq11}. Later we will present and discuss some notable 
applications of the theorem for special classes of nonlinearities of particular interest. (See 
Theorem \ref{LiouvilleThmEx} and Theorem \ref{LiouvilleThmEx-log} below.)

\begin{theorem} \label{thm Liouville}
Let $(M, g, e^{-f} dv_g)$ be a complete smooth metric measure space satisfying ${\mathscr Ric}_f (g) \ge 0$ 
on $M$. Let $u$ be a positive bounded solution to the equation $\Delta_f u + \Sigma(u) =0$. Suppose that 
along the solution $u$ we have $[1-\beta(\alpha/2+1)]\Sigma(u)-u\Sigma_u(u) \ge 0$ everywhere in $M$ 
for some $\alpha \ge 0$ and $0<\beta<1/(1+\alpha)$. Then $u$ 
must be a constant. In particular $\Sigma(u) =0$.
\end{theorem}

Note that since here the curvature bound is expressed as ${\mathscr Ric}_f(g) \ge 0$, the positivity of $u$ 
alone is not enough to grant the constancy of $u$ as can be seen from the 
discussion on $f$-harmonic functions above [\eqref{eq11} with $\Sigma=\Sigma(x,u) \equiv 0$]. For 
analogous results under the stronger curvature bound ${\mathscr Ric}_f^m(g) \ge 0$ but without the 
boundedness assumption on $u$ see \cite{TVahGrad}.

The proof of Theorem \ref{thm Liouville} is a direct consequence of the gradient estimates established in 
Theorems \ref{thm18} and \ref{thm18-global} and is presented in Section \ref{Liouville-Section}. We end 
this section by giving some nice applications of the above theorem. Towards this end consider first a 
superposition of power-like nonlinearities with constant coefficients $\mathsf{A}_j, \mathsf{B}_j$ and 
real exponents $p_j, q_j$ for $1 \le j \le N$ in the form 
\begin{equation} \label{polies PR}
\Sigma(u) = \mathscr{P}(u) + \mathscr{R}(u) 
= \sum_{j=1}^N \mathsf{A}_j u^{p_j} + \sum_{j=1}^N \mathsf{B}_j u^{q_j}.
\end{equation}
A direct calculation for the quantity pertaining to $\mathsf{T}_\Sigma(u)$ in Theorem \ref{thm Liouville} gives  
\begin{align}
[1-\beta(\alpha/2+1)]\Sigma(u)-u\Sigma_u(u) 
=&~\sum_{j=1}^N \mathsf{A}_j [1-\beta(\alpha/2+1) - p_j] u^{p_j} \nonumber \\
&+ \sum_{j=1}^N \mathsf{B}_j [1-\beta(\alpha/2+1) - q_j] u^{q_j} 
\end{align}
which is then easily seen to be non-negative, as required by Theorem \ref{thm Liouville}, upon suitably restricting the 
ranges of $\mathsf{A}_j$, $\mathsf{B}_j$ and $p_j, q_j$ as formulated below. This conclusion improves and extends 
earlier results on Yamabe type problems ({\it cf.} \cite{Dung, GidSp, Yang, Wu18}). Further applications and results 
in this direction will be discussed in a forthcoming paper ({\it see} also \cite{Taheri-GE-1, Taheri-GE-2}).

\begin{theorem} \label{LiouvilleThmEx}
Let $(M, g, e^{-f}dv_g)$ be a complete smooth metric measure space with 
${\mathscr Ric}_f(g) \ge 0$. Let $u$ be a positive bounded solution to the nonlinear elliptic equation   
\begin{equation} \label{elliptic PDE 2}
\Delta_f u + \sum_{j=1}^N \mathsf{A}_j u^{p_j} + \sum_{j=1}^N \mathsf{B}_j u^{q_j} = 0.
\end{equation}
Assume $\mathsf{A}_j \ge 0$, $\mathsf{B}_j \le 0$ and $p_j \le 1-\beta(\alpha/2+1)$, $q_j \ge 1-\beta(\alpha/2+1)$ 
for $1 \le j \le N$ and for some $\alpha \ge 0$ and $0<\beta < 1/(1+\alpha)$. Then $u$ must be a constant and $\Sigma(u)=0$. 
\end{theorem}

\begin{remark}
Note that the upper and lower bounds on the exponents $p_j, q_j$ are given by the same quantity 
$1-\beta(\alpha/2+1)$ which can be adjusted by optimising the parameters $\alpha, \beta$ within 
their respective range. In particular if $\mathscr{R}(u) \equiv 0$ we only require $\mathsf{A}_j \ge 0$ 
and $p_j <1$ (with $\alpha=0$, $\beta \searrow 0$) and if $\mathscr{P}(u) \equiv 0$ we only require 
$\mathsf{B}_j \le 0$ and $q_j >0$ (with $\alpha=0$, $\beta \nearrow 1$). 
\end{remark}

As another application consider a superposition of logarithmic type and power-like nonlinearities 
with real exponents $p$, $q$, $s$ and constant coefficients $\mathsf{A}$, $\mathsf{B}$ and 
$\gamma \in \mathscr{C}^1(\mathbb R)$ in the form 
\begin{equation}
\Sigma(u) = u^s \gamma(\log u) + \mathsf{A} u^p + \mathsf{B} u^q. 
\end{equation}
A straightforward calculation for the quantity pertaining to $\mathsf{T}_\Sigma(u)$ in Theorem 
\ref{thm Liouville} then gives 
\begin{align}
[1-\beta(\alpha/2+1)] \Sigma (u) - u \Sigma_u (u)
&= u^s ([1-\beta(\alpha/2+1) -s] \gamma - \gamma') \\
&+ \mathsf{A} [1-\beta(\alpha/2+1)-p] u^p 
+ \mathsf{B} [1-\beta(\alpha/2+1) - q] u^{q}, \nonumber 
\end{align}  
and so a discussion similar to that given before along with an application of Theorem \ref{thm Liouville} 
leads to the following result.

\begin{theorem} \label{LiouvilleThmEx-log}
Let $(M, g, e^{-f}dv_g)$ be a complete smooth metric measure space with ${\mathscr Ric}_f(g) \ge 0$. 
Let $u$ be a positive bounded solution to the nonlinear elliptic equation   
\begin{equation} \label{elliptic PDE 3}
\Delta_f u + u^s \gamma(\log u) + \mathsf{A} u^p + \mathsf{B} u^q = 0.
\end{equation}
Assume that along the solution $u$ we have the inequality $\gamma'+[\beta(\alpha/2+1)+s-1] \gamma \le 0$ 
along with $\mathsf{A} \ge 0$, $\mathsf{B} \le 0$, $p \le 1-\beta(\alpha/2+1)$ and $q \ge 1-\beta(\alpha/2+1)$ 
for some $\alpha \ge 0$ and $0<\beta < 1/(1+\alpha)$. Then $u$ must be a constant and $\Sigma(u)=0$. 
\end{theorem}

We end this section with the following global Hamilton-type result for equation \eqref{eq11} with a dimension 
free constant. With the aid of this bound one can then prove a global Harnack-type inequality. Note that here 
$M$ is assumed to be closed.

\begin{theorem} \label{global-Hamilton-Thm}
Let $(M, g, e^{-f} dv_g)$ be a complete closed smooth metric measure space satisfying the global curvature lower 
bound ${\mathscr Ric}_f (g) \ge - \mathsf{k} g$ with $\mathsf{k} \ge 0$. Suppose u is a bounded positive solution 
to $\Delta_f u +\Sigma(u)=0$ with $0<u\le A$. Moreover assume that along the solution $u$ we have 
$\Sigma(u) \ge 0$ and $\Sigma(u)/u - 2 \Sigma'(u) \ge 0$. Then 
\begin{align} \label{eq7.11}
|\nabla \log u|^2  \le 2 \mathsf{k} [1+ \log(A/u)]. 
\end{align}
Furthermore, there holds the  following interpolation-type Harnack inequality 
\begin{align}\label{eq16}
0<u(p) \le (eA)^{\varepsilon/(1+\varepsilon)} e^{\mathsf{k} d^2(p,q)/(2 \varepsilon)} 
[u(q)]^{1/(1+\varepsilon)}, 
\end{align}
for $\varepsilon >0$, $p, q \in M$.
\end {theorem}

\section{Proof of the Souplet-Zhang estimate in Theorem \ref{thm1}}\label{sec3}

Before proceeding to the proof of Theorem \ref{thm1} we need some intermediate 
results. Lemma \ref{lem20} and Lemma \ref{lem21} below ultimately lead to an elliptic differential 
inequality for a suitable quantity built out of the solution $u$ in Lemma \ref{lem2112} that plays a key role 
in the proof of the theorem.

\begin{lemma} \label{lem20}
Let $u$ be a positive bounded solution to \eqref{eq11} with $0<u \le A$. Then the function 
$h = \log (u/A)$ satisfies the equation 
\begin{equation} \label{h evolution eq}
\Delta_f h + |\nabla h|^2 + A^{-1} e^{-h} \Sigma(x,A e^h) =0.
\end{equation} 
\end{lemma}

\begin{proof} 
This is an easy calculation and the proof is left to the reader. 
\end{proof}

\begin{lemma} \label{lem21}
Let $u$ be a positive solution to \eqref{eq11} with $0<u \le A$. 
Put $h = \log (u/A)$ and let $H = |\nabla h|^2/(1-h)^2$. Then $H$ satisfies the equation  
\begin{align}\label{eq21}
\Delta_f H =&~\frac{2 {\mathscr Ric}_f}{(1-h)^2} (\nabla h, \nabla h) 
+ \frac{2h \langle \nabla h, \nabla H \rangle}{1-h} 
+ 2 (1-h) H^2 \nonumber \\
&+ 2 \left| \frac{\nabla^2 h}{1-h} + \frac{\nabla h \otimes \nabla h}{(1-h)^2} \right|^2 
- \frac{2 \langle \nabla h, \Sigma_x(x,Ae^h) \rangle}{A e^h (1-h)^2} \nonumber \\
&- 2 H \left( \Sigma_u (x,Ae^h) + \frac{h \Sigma (x,Ae^h)}{Ae^h(1-h)} \right).
\end{align}
\end{lemma}

\begin{proof}
In view of the relation $\nabla \log (1-h) = -\nabla h /(1-h)$ we have  
$H= |\nabla \log (1-h)|^2$. Moreover an easy calculation gives 
$\nabla^2 \log (1-h) = - \nabla^2h/(1-h) - (\nabla h \otimes \nabla h)/(1-h)^2$ 
and $-\Delta_f \log (1-h) = \Delta_f h / (1-h) + |\nabla h|^2/(1-h)^2$. Now an 
application of the Bochner-Weitzenb\"ock identity gives 
\begin{align} \label{eq3.4} 
\Delta_f H =&~\Delta_f|\nabla \log (1-h)|^2 \nonumber \\
=&~2|\nabla^2 \log (1-h)|^2 + 2 \langle \nabla \log (1-h) , \nabla \Delta_f \log (1-h) \rangle \nonumber \\
&+ 2{\mathscr Ric}_f(\nabla \log (1-h), \nabla \log (1-h)) \nonumber\\
=&~ \frac{2|\nabla^2 h|^2}{(1-h)^2}+ \frac{2 |\nabla h|^4}{(1-h)^4}+  \frac{2\langle \nabla h , \nabla |\nabla h|^2 \rangle}{(1-h)^3}\nonumber\\
& + 2 \left\langle \frac{\nabla h}{1-h} , \nabla \left[\frac{\Delta_f h }{1-h}+ \frac{|\nabla h|^2}{(1-h)^2}\right]\right \rangle
+ \frac{2{\mathscr Ric}_f(\nabla h, \nabla h)}{(1-h)^2}.
\end{align}

Next referring to the first term in the expression on the last line in \eqref{eq3.4} we can write 
\begin{align}\label{eq3.5}
\bigg\langle \frac{\nabla h}{1-h}, \nabla \bigg[\frac{\Delta_f h }{1-h} +& \frac{|\nabla h|^2}{(1-h)^2}\bigg]\bigg \rangle \nonumber\\
=&~\left \langle \frac{\nabla h}{1-h}, \frac{\nabla \Delta_f h}{1-h} + \frac{\Delta_f h \nabla h  }{(1-h)^2} \right \rangle \nonumber\\
&+ \left \langle \frac{\nabla h}{1-h}, \frac{\nabla |\nabla h|^2}{(1-h)^2} + \frac{2 |\nabla h|^2 \nabla h}{(1-h)^3}\right \rangle \nonumber\\
= &~ \frac{\langle \nabla h , \nabla \Delta_f h \rangle}{(1-h)^2} + \frac{ |\nabla h|^2 \Delta_f h}{(1-h)^3} 
+ \frac{\langle \nabla h , \nabla |\nabla h|^2 \rangle}{(1-h)^3} + \frac{2 |\nabla h|^4}{(1-h)^4}.
\end{align}
As a result by substituting \eqref{eq3.5} back in \eqref{eq3.4} it follows that  
\begin{align}\label{eq23}
\Delta_f H =&~\frac{2|\nabla^2 h|^2}{(1-h)^2}  + \frac{2\langle \nabla h , \nabla \Delta_f h \rangle}{(1-h)^2} 
+ \frac{2{\mathscr Ric}_f(\nabla h, \nabla h)}{(1-h)^2} \nonumber\\
& + \frac{4 \langle \nabla h , \nabla |\nabla h|^2 \rangle}{(1-h)^3}+ \frac{2 |\nabla h|^2 \Delta_f h}{(1-h)^3} 
+ \frac{6 |\nabla h|^4}{(1-h)^4}.
\end{align}

Now as for the term $\Delta_f h$ appearing on the right in \eqref{eq23} by substitution using the equation 
$\Delta_f h + |\nabla h|^2 + A^{-1} e^{-h} \Sigma=0$ it follows that  
\begin{align}
\Delta_f H  
=&~\frac{2 |\nabla^2 h|^2}{(1-h)^2} 
- \frac{2 \langle \nabla h, \nabla |\nabla h|^2 \rangle}{(1-h)^2} + \frac{2 \Sigma |\nabla h|^2}{Ae^h (1-h)^2} \nonumber \\
& - \frac{2 \langle \nabla h , \nabla \Sigma \rangle}{Ae^h (1-h)^2} 
+ \frac{2 {\mathscr Ric}_f (\nabla h, \nabla h)}{(1-h)^2}  + \frac{4\langle \nabla h , \nabla |\nabla h|^2 \rangle}{(1-h)^3} \nonumber \\
& - \frac{2 |\nabla h|^4}{(1-h)^3} -\frac{2 \Sigma |\nabla h|^2}{A e^h (1-h)^3} + \frac{6 |\nabla h|^4}{(1-h)^4}, 
\end{align}
and therefore a rearrangement of terms and basic considerations including the identity $\nabla \Sigma = \Sigma_x + Ae^h \Sigma_u \nabla h$, leads to 
\begin{align}\label {eq3.8}
\Delta_f H 
=&~\frac{ 2 {\mathscr Ric}_f(\nabla h, \nabla h)}{(1-h)^2} 
+ 2 \left| \frac{\nabla^2 h}{1-h} + \frac{\nabla h \otimes \nabla h}{(1-h)^2} \right|^2 + \frac{2 |\nabla h|^4}{(1-h)^3} \nonumber \\
&- \frac{2 \langle \nabla h, \nabla |\nabla h|^2 \rangle}{(1-h)^2} 
- \frac{4 |\nabla h|^4}{(1-h)^3}+ \frac{2 \langle \nabla h, \nabla |\nabla h|^2 \rangle}{(1-h)^3} + \frac{4 |\nabla h|^4}{(1-h)^4} \nonumber \\
& + \frac{2 \Sigma |\nabla h|^2}{Ae^h (1-h)^2}-\frac{2 \Sigma |\nabla h|^2}{A e^h (1-h)^3} 
-\frac{2 \langle \nabla h , \Sigma_x \rangle}{Ae^h(1-h)^2} -\frac{2 \Sigma_u |\nabla h|^2}{(1-h)^2}.
\end{align}
Finally, making note of the identity $(1-h)^3 \langle \nabla h, \nabla H \rangle=  (1-h) \langle \nabla h, \nabla |\nabla h|^2 \rangle + 2 |\nabla h|^4$,
the second line in \eqref{eq3.8} reduces to $2h\langle \nabla h, \nabla H \rangle / (1-h)$. Likewise recalling the relation 
$H=|\nabla h|^2/(1-h)^2$ and substituting back gives   
\begin{align}
\Delta_f H =& \, \frac{ 2 {\mathscr Ric}_f(\nabla h, \nabla h)}{(1-h)^2} 
+ 2 \left| \frac{\nabla^2 h}{1-h} + \frac{\nabla h \otimes \nabla h}{(1-h)^2} \right|^2  \nonumber \\
& +\frac{2h}{1-h} \langle \nabla h , \nabla H \rangle + 2(1-h) H^2 \nonumber\\
& -2H \left[ \Sigma_u + \frac{h\Sigma}{Ae^h (1-h)} \right] - \frac{2 \langle \nabla h ,\Sigma_x \rangle}{A e^h (1-h)^2}, 
\end {align}
which is the desired conclusion. The proof is thus complete. 
\end{proof}

\begin{lemma}\label{lem2112}
Under the assumptions of Lemma $\ref{lem21}$ and the curvature lower bound ${\mathscr Ric}_f(g) \ge -(n-1)kg $ with $k\ge 0$, 
the function $H = |\nabla h|^2/(1-h)^2$ satisfies the elliptic differential inequality 
\begin{align}\label{eq2112}
\Delta_f  H 
\ge &~2(1-h)H^2 - 2 (n-1)kH + \frac{2h \langle\nabla h, \nabla H\rangle}{1-h} \nonumber \\
&- \frac{2 \langle \nabla h, \Sigma_x(x,Ae^h) \rangle}{A e^h (1-h)^2} 
 - 2 H \left[\Sigma_u (x,Ae^h) + \frac{h \Sigma (x,Ae^h)}{Ae^h(1-h)} \right].
\end{align}
\end{lemma}

\begin{proof}
This is a straightforward consequence of the identity \eqref{eq21} in Lemma \ref{lem21} and the curvature lower bound 
${\mathscr Ric}_f(g) \ge -(n-1)kg$ in the lemma.
\end{proof}

\qquad \\
{\it Proof of Theorem $\ref{thm1}$.} 
In order to prove Theorem \ref{thm1}, the elliptic differential inequality in Lemma \ref{lem2112} will be combined with 
a localisation and cut-off function argument along with an application of the maximum principle to a suitably defined 
localised function. Now pick a reference point $p$ in $M$, fix $R \ge 2$ and let $r=r_p(x)$ denote the 
geodesic radial variable with respect to $p$. Below we use the radially symmetric cut-off function  
\begin{equation} \label{cut-off def}
\phi(x) = \bar{\phi} (r(x)),
\end{equation}
supported in $\mathcal{B}_{R} = \mathcal{B}_R(p) \subset M$. The function $\bar{\phi}=\bar \phi(r)$ appearing on the 
right-hand side of \eqref{cut-off def} is chosen to satisfy the following properties 
\begin{enumerate}[label=$(\roman*)$]
\item $\bar{\phi} \in \mathscr{C}^2([0,\infty), \mathbb{R})$, ${\rm supp} \, \bar{\phi} \subset [0,R]$ 
and $0 \le \bar{\phi} \le 1$, 
\item $\bar{\phi} \equiv 1$ in $[0,R/2] $ and $\bar{\phi}^{'} \equiv 0$ in $[0,R/2] \cup [R, \infty)$, 
\item for every $0<\varepsilon<1$ there exists $c_\varepsilon>0$ such that the bounds,  
\begin{equation}
-c_\varepsilon \bar{\phi}^\varepsilon/R \le \bar{\phi}' \leq 0, \qquad |\bar{\phi}''| \leq c_\varepsilon \bar{\phi}^\varepsilon / R^2, 
\end{equation} 
hold on $[0, \infty)$.

\end{enumerate}

Our aim is now to establish the desired estimate at every $x$ in ${\mathcal B}_{R/2}$ and for this we consider 
the localised function $\phi H$ with $\phi$ as in \eqref{cut-off def} and $H$ as in Lemma \ref{lem21}. To this end we 
start with $\Delta_f (\phi H) = \phi \Delta_f  H +  2\langle \nabla H , \nabla \phi \rangle+ H \Delta_f  \phi$ or after a 
slight adjustment, 
\begin{equation}
\Delta_f (\phi H) = \phi \Delta_f H +   2 [\langle\nabla \phi ,\nabla (\phi H) \rangle- |\nabla \phi|^2 H]/\phi  + H \Delta_f  \phi.
\end{equation}  
An application of the inequality \eqref{eq2112} in Lemma \ref{lem2112} and another adjustment of terms as above then gives 
\begin{align} \label{eq25}
\Delta_f (\phi H) 
\ge & \left\langle \left[ \frac{2h\nabla h}{1 -h} + \frac{2\nabla \phi}{\phi} \right],  \nabla (\phi H) \right\rangle
- \left\langle H \left[ \frac{2h \nabla h}{1 -h} + \frac{2 \nabla \phi}{\phi} \right], \nabla \phi \right\rangle \nonumber \\
&+ 2(1-h) \phi H^2 + H [\Delta_f -2(n-1) k] \phi 
- \frac{2 \phi \langle \nabla h, \Sigma_x(x,Ae^h) \rangle}{A e^h (1-h)^2} \nonumber \\
& - 2  \phi H \left[ \Sigma_u (x,Ae^h) + \frac{h \Sigma (x,Ae^h)}{Ae^h(1-h)} \right]. 
\end{align}

Assume now that $\phi H$ is maximised on $\mathcal{B}_{R}$ at the point $q$. Since 
the function $r=r(x)$ is only Lipschitz continuous at the cut locus of $p$, using a standard 
argument of Calabi ({\it see}, e.g., \cite{SchYau} p.~21), we can assume without loss of 
generality that $q$ is not in the cut locus of $p$ and hence $\phi H$ is smooth at $q$ 
for the application of the maximum principle. Moreover we assume $(\phi H)(q) >0$ 
as otherwise the result is trivial in view of $H(x) \le 0$ in $\mathcal{B}_{R/2}$. Now 
at the maximal point $q$  we have $\Delta_f(\phi H) \le 0$ and $\nabla(\phi H) =0$. 
From \eqref{eq25} it thus follows that  
\begin{align} 
2(1-h)\phi H^2 
\le& \left\langle H \left[ \frac{2h \nabla h}{1 -h} + \frac{2 \nabla \phi}{\phi} \right], \nabla \phi \right\rangle 
- H \Delta_f \phi + 2(n-1) kH\phi \nonumber \\
&+ \frac{2 \phi \langle \nabla h, \Sigma_x(x,Ae^h) \rangle}{A e^h (1-h)^2} 
+ 2 \phi H \left[ \Sigma_u (x,Ae^h) + \frac{h \Sigma (x,Ae^h)}{Ae^h(1-h)} \right],
\end{align}
at $q$ or dividing through by $2(1-h) \ge 0$ that 
\begin{align}\label{eq27}
\phi H^2 
\le & \left\langle 
H \left[ \frac{h \nabla h}{1 -h} + \frac{\nabla \phi}{\phi} \right], \frac{\nabla \phi}{1-h} \right\rangle 
+ \frac{H[-\Delta_f \phi + 2 (n-1)k\phi]}{2(1-h)} \nonumber \\
&+ \frac{\phi \langle \nabla h, \Sigma_x(x,Ae^h) \rangle}{A e^h (1-h)^3} 
+ \phi H \left[ \frac{\Sigma_u (x,Ae^h)}{1-h} + \frac{h \Sigma (x,Ae^h)}{Ae^h(1-h)^2} \right].
\end{align}
The goal is now to use \eqref{eq27} to establish the required estimate at $x$. To this end we 
consider two separate cases. Firstly, the case $d(q) \le 1$ and next the case 
$d(q) \ge 1$. \\
{\bf Case 1.} Since $\phi \equiv 1$ in ${\mathcal B}_{R/2}$ [i.e., for all $x$ with $d(x) \le R/2$, 
where $R \ge 2$ by property $(ii)$] all the terms involving derivatives of $\phi$ at $q$ vanish (in particular $\nabla \phi=0$ and
$\Delta_f \phi =0$). So as a result it follows from (\ref{eq27}) that at 
the point $q$, we have the bound
\begin{align*}
\phi H^2 
\le \frac{(n-1)k\phi H}{(1-h)} + \frac{\phi |\langle \nabla h, \Sigma_x(x,Ae^h) \rangle|}{A e^h (1-h)^3} 
+ \phi H \left[ \frac{\Sigma_u (x,Ae^h)}{1-h} + \frac{h \Sigma (x,Ae^h)}{Ae^h(1-h)^2} \right]_+,
\end{align*}
and so 
\begin{align*}
\phi H^2 \le \frac{(n-1)k \sqrt \phi H}{(1-h)} 
+ \phi \sqrt H \frac{|\Sigma_x(x,Ae^h)|}{A e^h (1-h)^2}
+ \sqrt \phi H \left[ \frac{\Sigma_u (x,Ae^h)}{1-h} + \frac{h \Sigma (x,Ae^h)}{Ae^h(1-h)^2} \right]_+. \nonumber
\end{align*}
By an application of Young's inequality it then follows after rearranging terms that
\begin{align*}
\phi H^2 \le C \left\{ k^2 
+ \left[ \frac{\Sigma_u (x,Ae^h)}{1-h} + \frac{h \Sigma (x,Ae^h)}{Ae^h(1-h)^2} \right]^2_+
+ \left[ \frac{|\Sigma_x(x,Ae^h)|}{A e^h (1-h)^2} \right]^{4/3} \right\}.
\end{align*}
As $\phi \equiv 1$ when $d(x) \le R/2$, we have $H(x)=[\phi H](x) \le [\phi H](q) \le [\sqrt \phi H] (q)$. 
Hence by recalling $H=|\nabla h|^2/(1-h)^2$ and $h=\log(u/A)$, we arrive at the bound at $x$ 
\begin{align}
\frac{\left| \nabla h \right|}{1- h} \le C \left\{ \sqrt k 
+ \sup_{\mathcal{B}_{R}} \left[ \frac{Ae^h (1-h) \Sigma_u + h \Sigma}{Ae^h(1-h)^2} \right]_+^{1/2} 
+ \sup_{\mathcal{B}_{R}} \left[ \frac{|\Sigma_x|}{Ae^h (1-h)^2} \right]^{1/3} \right\}.   
\end{align}
Evidently this a special case of the inequality \eqref{eq13} upon noting \eqref{eq2.2} and \eqref{eq2.3} 
and so the proof of the estimate is complete in this case.

\qquad \\
{\bf Case 2.} Upon referring to the right-hand side of \eqref{eq27}, and noting the properties of $\bar \phi$ as listed at the 
start we proceed onto bounding the full expression on the right-hand side on \eqref{eq27} in the case $d(q) \ge 1$. 
Towards this end dealing with the first term first, we have 
\begin{align} \label{eq28}
H \left\langle \left[ \frac{h \nabla h}{1 -h} + \frac{\nabla \phi}{\phi} \right], \frac{\nabla \phi}{1-h} \right\rangle
&\le H \left[ \frac{h |\nabla h|}{1 -h} + \frac{|\nabla \phi|}{\phi} \right] \frac{|\nabla \phi|}{1-h} \nonumber \\
&\le H \left[ h \sqrt H + \frac{|\nabla \phi|}{\phi} \right] \frac{|\nabla \phi|}{1-h} \\
&\le H \left [ \phi^{1/4} \frac{\sqrt H |h|}{1-h} \frac{|\nabla \phi|}{\phi^{3/4}} + \frac{|\nabla \phi|^2}{\phi^{3/2}} \right] \sqrt \phi 
\le \frac{\phi H^2}{4} + \frac{C}{R^4}. \nonumber 
\end{align}
In much the same way regarding the terms involving $\Sigma$ upon setting $\mathsf{P}_\Sigma(u)=|\Sigma_x|/(Ae^h)$ we 
have firstly 
\begin{align}
\frac{\phi \langle \nabla h, \Sigma_x \rangle}{A e^h (1-h)^3} \le \frac{\phi |\nabla h| |\Sigma_x|}{A e^h (1-h)^3} 
&= \frac{\phi \sqrt H |\Sigma_x|}{A e^h (1-h)^2} \nonumber \\
&\le \frac{\phi H^2}{8} + C \left(\frac{|\Sigma_x|}{Ae^h} \right)^{4/3}
= \frac{\phi H^2}{8} + C \mathsf{P}_\Sigma^{4/3} (u), 
\end{align}
and likewise for the subsequent terms, upon noting $-1 \le h/(1-h) \le 0$, $h \le 0$ and $0 \le \phi \le 1$, we have
\begin{align}
\phi H \left[ \frac{\Sigma_u}{1-h} + \frac{h \Sigma}{Ae^h(1-h)^2} \right] &= \phi H \left[ \frac{Ae^h (1-h) \Sigma_u 
+ h \Sigma}{Ae^h(1-h)^2} \right] \\
&\le \frac{\phi H^2}{8} + C \left[ \frac{Ae^h (1-h) \Sigma_u + h \Sigma}{A e^h (1-h)^2} \right]^2_+ = \frac{\phi H^2}{8} 
+ C \mathsf{R}^2_\Sigma (u), \nonumber 
\end{align}
where in the last equation we have set 
\begin{equation}
\mathsf{R}_\Sigma (u) = \{ [Ae^h (1-h) \Sigma_u + h \Sigma]/[(1-h)^2 Ae^h] \}_+.
\end{equation}

Now for the term $\Delta_f \phi$ we use the Wei-Wylie weighted Laplacian comparison theorem taking advantage 
of the fact that it only depends on the lower bound on ${\mathscr Ric}_f(g)$ (\cite{[WeW09]}). Indeed recalling 
$1 \le d(q) \le R$, it follows from ${\mathscr Ric}_f(g) \ge - (n-1) k g$ with $k \ge 0$ and Theorem 3.1 in 
\cite{[WeW09]} that [{\it see} \eqref{sigma def alt} for notation]
\begin{equation} \label{fLap-of-r}
\Delta_f r (x) \le \gamma_{\Delta_f} +(n-1) k (R-1),
\end{equation}
whenever $1 \le r \le R$ [in particular at the point $q$]. 
Thus proceeding on to bounding $-\Delta_f \phi$, upon referring to \eqref{cut-off def} and using $(ii)$ [$\bar \phi^{'} =0$ 
when $0 \le r \le R/2$], $(iv)$ [$\bar \phi^{'} \le 0$ when $0 \le r < \infty$] 
we have $-\Delta_f \phi = -(\bar \phi^{''}|\nabla r|^2 + \bar \phi^{'} \Delta_f r)$, and so 
\begin{align} \label{bound Delta f alpha eq}
- \Delta_f \phi 
&\le \left[ \frac{|\bar \phi^{''}|}{\sqrt{\bar \phi}} + \frac{|\bar \phi^{'}|}{\sqrt{\bar \phi}} 
([\gamma_{\Delta_f}]_+ + (n-1) k (R-1) ) \right] \sqrt{\bar \phi} \nonumber \\
&\le C \left(  \frac{1}{R^2} + \frac{[\gamma_{\Delta_f}]_+}{R} + k \right) \sqrt{\phi}.
\end{align}
Next referring to the last term on the first line in \eqref{eq27}, the above, along with $1-h \ge 1$, 
after an application of Young inequality gives
\begin{align} \label{eq28b} 
\frac{H [-\Delta_f \phi + 2(n-1)k \phi]}{2(1-h)} \le \frac{\phi H^2}{4} + C \left( \frac{1}{R^4} 
+ \frac{[\gamma_{\Delta_f}]_+^2}{R^2} + k^2 \right). 
\end{align}

Thus, by putting together the above fragments, namely, the bounds in \eqref{eq28}--\eqref{eq28b}, 
we obtain, after reverting to $u=Ae^h$, the following upper bound on $\phi H^2$ at $q$,  
\begin{align} \label{middle-final}
\phi H^2 \le&  \, C \left\{ \frac{1}{R^4}+\frac{[\gamma_{\Delta_f}]_+^2}{R^2} + k^2 
+ \sup_{\mathcal{B}_{R}} \left[ \mathsf{R}^{2}_\Sigma(u) + \mathsf{P}^{4/3}_\Sigma(u) \right] \right\}.
\end{align}

Finally, recalling that $\phi H$ is maximised on $\mathcal{B}_R$ at $q$, $\phi \equiv 1$ 
on $\mathcal{B}_{R/2}$ and $0\le\phi\le1$ on $\mathcal{B}_R$ we can write 
\begin{equation}
\sup_{\mathcal{B}_{R/2}} H^2 = \sup_{\mathcal{B}_{R/2}} (\phi^2 H^2) 
\le \sup_{\mathcal{B}_{R}} (\phi^2 H^2)  = (\phi^2 H^2)(q) \le (\phi H^2)(q).  
\end{equation} 
Hence upon noting $H = |\nabla h|^2/(1-h)^2$ and $h=\log u$, the above result in 
\begin{align} \label{eq214}
\frac{\left| \nabla \log u \right|}{1- \log(u/A)}
\le& C \left\{ \frac{1}{R}+\sqrt{\frac{[\gamma_{\Delta_f}]_+}{R}} + \sqrt{k} 
+ \sup_{\mathcal{B}_{R}} \left[ \mathsf{R}^{1/2}_\Sigma(u) + \mathsf{P}^{1/3}_\Sigma(u) \right]
 \right\}.
\end{align}
Thus in either of the two cases above we have shown the estimate to be valid at $x$. The arbitrariness of 
$x$ in ${\mathcal B}_{R/2}$ gives the required conclusion and completes the proof. \hfill $\square$

\begin{remark} \label{Ricmf-proof-1-remark} {
Under the stronger curvature assumption ${\mathscr Ric}_f^m(g) \ge -(m-1)kg$ with $n \le m < \infty$ 
and $k \ge 0$, we have from the Wei-Wylie weighted Laplacian comparison theorem  
({\it see} also Corollary 1.2 in \cite{LiX} as well as \cite{Lott, Qian1}) the $f$-Laplacian inequality 
\begin{equation}
\Delta_f r \le (m-1) \sqrt k \coth (\sqrt k r).
\end{equation} 
[Compare with \eqref{fLap-of-r}.]  
Therefore, for the cut-off function $\phi$ in \eqref{cut-off def}, by virtue of $\phi$ being radial and 
$\bar \phi' \le 0$, we can deduce the lower bound 
\begin{equation*}
\Delta_f \phi = \bar \phi'' |\nabla r|^2 + \bar \phi' \Delta_f r 
\ge \bar \phi'' + (m-1) \bar \phi' \sqrt k \coth (\sqrt k r). 
\end{equation*}
From the bound $\sqrt k \coth (\sqrt k r) \le \sqrt k \coth (\sqrt k R/2) \le (2+\sqrt k R)/R$ 
for $R/2 \le r \le R$ (here we are making use of $v \coth v \le 1+v$ and the monotonicity of $\coth v$ for 
$v>0$) upon noting $\bar \phi' \equiv 0$ for $0 \le r \le R/2$ it follows that 
\begin{align*}
- \Delta_f \phi &\le -[\bar \phi'' + (m-1) \bar \phi' \sqrt k \coth (\sqrt k R/2)] 
\le |\bar \phi''| + (m-1)(2/R+\sqrt k) |\bar \phi'|. 
\end{align*}
Thus $- \Delta_f \phi \le c/R^2 + c(m-1)[2+\sqrt kR]/R^2 \le c[1+2(m-1)]/R^2+c(m-1) \sqrt k$ where in the 
second inequality we have made use of $R \ge 2$ [compare with \eqref{bound Delta f alpha eq}]. Making 
use of this inequality in \eqref{eq28b} instead of the previously used \eqref{bound Delta f alpha eq} implies 
that in \eqref{middle-final} and hence \eqref{eq214} we can remove the term \eqref{extra-term-Ricmf} at the cost of 
adjusting the constant $C>0$ (so that now it will depend on $m$ too).

}
\end{remark}

\section{Proof of the Harnack inequality in Theorem \ref{cor Harnack}} \label{sec4}

In this section as the title suggests we give a proof of the local and global Harnack inequalities in Theorem \ref{cor Harnack}. 
In its local form this is a consequence of the local estimate in Theorem \ref{thm1} and in its global form this is a consequence 
of the global estimate in Theorem \ref{thm1-global}.

Towards this end pick $x, y$ in $M$ and let $\zeta=\zeta(s)$ with $0 \le s \le 1$ be a shortest geodesic curve with 
respect to the metric $g$ joining $x, y$ in $\mathcal{B}_{R/2} \subset M$. Thus $\zeta(0)=x$, $\zeta(1)=y$ 
and $\zeta(s) \in \mathcal{B}_{R/2}$ for all $0 \le s \le 1$. 

By the compactness of ${\mathcal B}_{R} \subset M$ and $u>0$ 
we have $\sup_{{\mathcal B}_R} [\mathsf{R}_\Sigma(u)^{1/2} + \mathsf{P}_\Sigma(u)^{1/3}] < \infty$, 
and so in \eqref{exponent-local} we have $\gamma \in (0,1)$. Now invoking the local estimate 
in Theorem \ref{thm1} (recalling that $1-h \ge 1$) and with $d=d(x, y)$ we can write 
\begin{align}
\log \frac{1-h(y)}{1-h(x)} &= \int_0^1 \frac{d}{ds} \log [1-h(\zeta(s))] \, ds 
= \int_0^1 - \frac{\langle \nabla h (\zeta(s)), \zeta'(s) \rangle}{1-h(\zeta(s))} \, ds \nonumber \\
&\le \int_0^1 \frac{|\nabla h (\zeta(s))||\zeta'|}{1-h(\zeta(s))} \, ds 
\le \sup_{\mathcal{B}_{R/2}} \left[ \frac{|\nabla h|}{1-h} \right] \int_0^1 |\zeta'| \, ds \\
&\le C d \left\{ 
\frac{1}{R}+\sqrt{\frac{[\gamma_{\Delta_f}]_+}{R}} + \sqrt{k} 
+ \sup_{\mathcal{B}_{R}} \left[ \mathsf{R}^{1/2}_\Sigma(u) + \mathsf{P}^{1/3}_\Sigma(u) \right]
\right\}, \nonumber
\end{align}
where $\zeta'=d \zeta / ds$ and we have substituted $|\nabla h|/(1-h) = |\nabla \log (u/A)|/[1-\log(u/A)]$. 
Therefore a straightforward calculation gives, 
\begin{align}
\frac{\log [eA/u(y)]}{\log [eA/u(x)]} &= \frac{1- \log[u(y)/A]}{1- \log[u(x)/A]} = \frac{1-h(y)}{1-h(x)} \\
&\le {\rm exp} \left[ C d \left( 
\frac{1}{R}+\sqrt{\frac{[\gamma_{\Delta_f}]_+}{R}} + \sqrt{k} 
+ \sup_{\mathcal{B}_{R}} \left[ \mathsf{R}^{1/2}_\Sigma(u) + \mathsf{P}^{1/3}_\Sigma(u) \right]
\right) \right] \nonumber 
\end{align}
and so the assertion follows. \hfill $\square$

\section{Proof of the Hamilton estimate in Theorem \ref{thm18}}\label{sec5}

Before moving on to the proof of Theorem \ref{thm18} we pause to state and prove two useful lemmas that are needed 
later on in establishing an elliptic differential inequality for a suitable quantity built out of the solution $u$ [{\it see} 
\eqref{Pf thm18 eq1} in the proof of Theorem \ref{thm18}].

\begin{lemma} \label{lemma h one}
Let $u$ be a positive solution to the nonlinear elliptic equation \eqref{eq11} and put $h=u^\beta$ where $\beta \in (0, 1)$. 
Then $h$ satisfies the equation 
\begin{align} \label{eq4.1}
\Delta_f h + (1-\beta) |\nabla h|^2/(\beta h) 
+ \beta h^{1-1/\beta} \Sigma(x, h^{1/\beta}) =0.
 \end{align}
\end{lemma}

\begin{proof}
A basic calculation gives $\Delta_f h = \beta (\beta-1) u^{\beta-2} |\nabla u|^2 + \beta u^{\beta -1} \Delta_f u$ 
and so \eqref{eq4.1} follows by substitution using $u=h^{1/\beta}$, $|\nabla u|^2 = h^{2(1/\beta-1)} |\nabla h|^2/\beta^2$ and \eqref{eq11}. 
\end{proof}

\begin{lemma} \label{lemma h two}
Under the assumptions of Lemma $\ref{lemma h one}$ on $u$ and $h=u^\beta$ with $0<\beta<1$, the function 
$G=G_\alpha^\beta= h^\alpha |\nabla h|^2$ with $\alpha \ge 0$ satisfies the equation 
\begin{align}  \label{eq4.5}
\Delta_f G_\alpha^\beta = &~ 2 h^{\alpha}|\nabla ^2 h|^2 + \alpha (\alpha -1) h^{\alpha -2} |\nabla h|^4 
+ 2h^{\alpha} {\mathscr Ric}_f (\nabla h, \nabla h)\nonumber\\
& + \alpha h^{\alpha -1} |\nabla h|^2 \Delta_f h  + 2 h^\alpha \langle \nabla h , \nabla \Delta_f h \rangle 
+ 2 \alpha h^{\alpha-1} \langle \nabla h , \nabla |\nabla h|^2 \rangle.
 \end{align}
 \end{lemma}

\begin{proof}
In order to calculate the action of $\Delta_f$ on $G$ we first write $G = |\nabla \bar h_\alpha|^2$ 
where $\bar h_\alpha = 2 h^{\alpha/2 +1}/(\alpha +2)$. Indeed here we have the relations 
$\nabla \bar h_\alpha = h^{\alpha/2} \nabla h$ and  
\begin{align}
\Delta_f \bar h_\alpha &= \Delta \bar h_\alpha - \langle \nabla f , \nabla \bar h_\alpha \rangle 
= \text{div} (h^{\alpha/2} \nabla h) - h^{\alpha/2}\langle \nabla f , \nabla h \rangle \nonumber\\
&= (\alpha/2) h^{\alpha/2 -1} |\nabla h|^2 + h^{\alpha/2} \Delta h - h^{\alpha/2}\langle \nabla f , \nabla h \rangle \nonumber\\
& = (\alpha/2) h^{\alpha/2 -1} |\nabla h|^2 + h^{\alpha/2} \Delta_f h.
\end{align}
Moreover a straightforward calculation gives the Hessian 
\begin{equation}
{\rm Hess} (\bar h_\alpha) = \nabla^2 \bar h_\alpha = h^{\alpha/2}[\nabla ^2 h + \alpha(\nabla h \otimes \nabla h)/(2h)].
\end{equation} 

Now applying the weighted Bochner-Weitzenb\"ock formula \eqref{Bochner} to $\bar h_\alpha$ and recalling the relation 
$G = |\nabla \bar h_\alpha|^2$ gives 
\begin{align}
\Delta_f G = &~ \Delta_f |\nabla \bar h_\alpha|^2= 2 |\nabla^2 \bar h_\alpha|^2 
+ 2 \langle \nabla \bar h_\alpha, \nabla \Delta_f \bar h_\alpha \rangle 
+2 {\mathscr Ric}_f (\nabla \bar h_\alpha , \nabla \bar h_\alpha)\nonumber\\
= &~ 2 h^{\alpha}|\nabla ^2 h|^2 + \alpha ^2/2 h^{\alpha -2} |\nabla h|^4 
+ \alpha h^{\alpha -1} \langle \nabla h, \nabla |\nabla h|^2 \rangle \nonumber\\
& + 2 h^{\alpha/2} \langle \nabla h , \nabla [h^{\alpha/2} \Delta_f h 
+ \alpha /2 h^{\alpha/2 -1} |\nabla h|^2] \rangle 
+ 2 h^{\alpha} {\mathscr Ric}_f (\nabla h, \nabla h )\nonumber\\
= & ~ 2 h^{\alpha}|\nabla ^2 h|^2 + \alpha ^2/2 h^{\alpha -2} |\nabla h|^4 
+ \alpha h^{\alpha -1} \langle \nabla h, \nabla |\nabla h|^2 \rangle \nonumber\\
& + 2 h^{\alpha} \langle \nabla h, \nabla \Delta_f h \rangle + \alpha h^{\alpha -1} |\nabla h|^2\Delta_f h
+ \alpha ( \alpha/2 -1) h^{\alpha -2} |\nabla h|^4 \nonumber\\
&+ \alpha h^{\alpha -1} \langle \nabla h, \nabla |\nabla h|^2 \rangle 
+ 2 h^{\alpha} {\mathscr Ric}_f (\nabla h, \nabla h ).
\end{align}
A rearrangement of terms now results in  
\begin{align}
\Delta_f G = 
&~ 2 h^{\alpha}|\nabla ^2 h|^2 + 2 \alpha h^{\alpha-1} \langle \nabla h , \nabla |\nabla h|^2 \rangle 
+ 2 h^\alpha \langle \nabla h , \nabla \Delta_f h \rangle \nonumber\\
& + \alpha h^{\alpha -1} |\nabla h|^2 \Delta_f h + \alpha (\alpha -1) h^{\alpha -2} |\nabla h|^4 
+ 2h^{\alpha} {\mathscr Ric}_f (\nabla h, \nabla h), 
\end{align}
which is immediately seen to be the required conclusion.
\end{proof}

\qquad \\
{\it Proof Theorem $\ref{thm18}$.} 
From the curvature lower bound ${\mathscr Ric}_f (g)\ge - (n-1) kg $ and the inequality $h>0$, upon substituting into the identity 
\eqref{eq4.5} in Lemma \ref{lemma h two} and making use of \eqref{eq4.1} in Lemma \ref{lemma h one}, it follows that 
\begin{align} \label{Pf thm18 eq1}
\Delta_f G \ge &~2h^{\alpha} |\nabla^2 h|^2 + 2 \alpha h^{\alpha-1} \langle \nabla h , \nabla |\nabla h|^2 \rangle\nonumber\\
&+ 2(\beta-1)/\beta h^{\alpha-1}  \langle \nabla h, \nabla |\nabla h|^2 \rangle \nonumber \\ 
&-2(\beta-1)/\beta h^{\alpha -2} |\nabla h|^4 \nonumber\\
&- 2 \beta h^{\alpha} \langle \nabla h, \nabla [h^{1-1/\beta} \Sigma(x, h^{1/\beta})] \rangle\nonumber\\
& -2(n-1)k h^{\alpha} |\nabla h|^2 + \alpha(\alpha -1) h^{\alpha -2} |\nabla h|^4\nonumber \\
& +\alpha(\beta -1)/\beta h^{\alpha-2} |\nabla h|^4 \nonumber\\
& - \alpha \beta h ^{\alpha-1/\beta} |\nabla h|^2 \Sigma(x,h^{1/\beta}).
\end{align}
Abbreviating hereafter the arguments of $\Sigma$ and its partial derivatives for convenience we can proceed by writing 
\begin{align} \label{Pf thm18 eq3}
\langle \nabla h, \nabla [h^{1-1/\beta} \Sigma] \rangle 
= (\beta-1)/(\beta h^{1/\beta}) \Sigma|\nabla h|^2  + h^{1-1/\beta} \langle \nabla h, \nabla \Sigma \rangle, 
\end{align}
where the last term on the right-hand side of \eqref{Pf thm18 eq3} can in turn be calculated as 
\begin{align} \label{Pf thm18 eq4}
\langle \nabla h, \nabla \Sigma \rangle &= \langle \nabla h, \Sigma_x + \nabla (h^{1/\beta}) \Sigma_u \rangle \nonumber \\
&= \langle \nabla h, \Sigma_x \rangle + \beta^{-1} h^{1/\beta -1} |\nabla h|^2 \Sigma_u \nonumber \\
&= \langle \nabla h, \Sigma_x \rangle + \beta^{-1} h^{1/\beta -2} G \Sigma_u. 
\end{align}
Substituting the descriptions of the inner products \eqref{Pf thm18 eq3} and \eqref {Pf thm18 eq4} back in \eqref{Pf thm18 eq1} 
then result in the inequality 
\begin{align} \label{eq4.9}
\Delta_f G \ge &~2h^{\alpha} |\nabla^2 h|^2 + 2 \alpha h^{\alpha-1} \langle \nabla h , \nabla |\nabla h|^2\rangle\nonumber\\
& -2(n-1)k h^{\alpha} |\nabla h|^2 \nonumber\\
&+ 2(\beta-1)/\beta h^{\alpha-1}  \langle \nabla h, \nabla |\nabla h|^2 \rangle \nonumber\\
&- [(2-\alpha^2) \beta -2+\alpha]/\beta h^{\alpha -2}|\nabla h|^4 \nonumber \\ 
&-[(2+\alpha)\beta -2] h^{\alpha -1/\beta} \Sigma |\nabla h|^2 \nonumber\\
&-2\beta h^{1+\alpha -1/\beta} \langle \nabla h , \Sigma_x \rangle 
-2 h^{\alpha} |\nabla h|^2 \Sigma_u.
\end{align}

Now using $  h^\alpha |\nabla^2 h|^2 + \alpha h^{\alpha -1}\langle \nabla h, \nabla |\nabla h|^2 \rangle
=h^{\alpha -1}[|\sqrt{h} \nabla^2 h + \alpha[\nabla h \otimes \nabla h]/\sqrt{h}|^2 -\alpha^2 |\nabla h|^4 /h] \ge -\alpha^2 h^{\alpha -2} |\nabla h|^4 $, 
and substituting accordingly in the first line on the right-hand side of \eqref{eq4.9}, we can rewrite the latter inequality as 
\begin{align} \label{Pf thm18 eq2}
\Delta_f G \ge & - 2(n-1)k h^{\alpha} |\nabla h|^2+2(\beta-1)/\beta h^{\alpha -1}\langle \nabla h, \nabla |\nabla h|^2 \rangle \nonumber\\
&- [(2+\alpha^2) \beta -2+\alpha)]/\beta  h^{\alpha -2}|\nabla h|^4 \nonumber \\
& -[(2+\alpha)\beta -2] h^{\alpha -1/\beta} \Sigma |\nabla h|^2 \nonumber\\
&-2\beta h^{1+\alpha -1/\beta} \langle \nabla h , \Sigma_x \rangle 
-2 h^{\alpha} |\nabla h|^2 \Sigma_u.
\end{align}

Let us denote by $\mathsf{Z}_\Sigma$ the sum of the last three terms on the right-hand side of \eqref{Pf thm18 eq2}, and note upon recalling 
$G=h^{\alpha}|\nabla h|^2$ and $u=h^{1/\beta}$ that 
\begin{align} \label{Z-Sigma}
\mathsf{Z}_\Sigma= \frac{[2-(2+\alpha)\beta] \Sigma - 2u \Sigma_u}{u} G 
- \frac{2\beta h^{1+\alpha} \langle \nabla h, \Sigma_x \rangle}{u}.
\end{align}

Since $\langle \nabla h, \nabla G \rangle = \langle \nabla h, \nabla (h^{\alpha} |\nabla h|^2) \rangle =\alpha h^{\alpha -1} |\nabla h|^4 
+ h^\alpha \langle \nabla h, \nabla |\nabla h|^2 \rangle$ 
by substituting these back in the inequality \eqref{Pf thm18 eq2} we can then write  
\begin{align}
\Delta_f G 
\ge & - 2(n-1)k h^{\alpha} |\nabla h|^2+2(\beta-1)/\beta h^{\alpha -1} \langle \nabla h, \nabla |\nabla h|^2 \rangle + \mathsf{Z}_\Sigma\nonumber\\
&- [(2+\alpha^2) \beta -2+\alpha ]/\beta  h^{\alpha -2}|\nabla h|^4 \nonumber \\
= &- 2(n-1)k h^{\alpha} |\nabla h|^2 +2(\beta-1)/(\beta h) \langle \nabla h, \nabla (h^{\alpha}|\nabla h|^2) \rangle + \mathsf{Z}_\Sigma  \nonumber \\
&+ [2 -(2+\alpha^2)\beta-(2\beta-1)\alpha ]/\beta  h^{\alpha -2}|\nabla h|^4 \nonumber \\
\ge & - 2(n-1)k G+ 2(\beta-1)/(\beta h) \langle \nabla h, \nabla G \rangle  + \mathsf{Z}_\Sigma\nonumber\\
&+ [2 -(2+\alpha^2)\beta-(2\beta-1)\alpha ]/(\beta  h^{\alpha +2}) G^2.
\end{align}

Next, localising by taking a space-time cut-off function $\phi$ as in \eqref{cut-off def} and following a similar procedure to those 
used in the proof of Theorem \ref{thm1}, we can write 
\begin{align}  \label{Pf thm18 eq5}
\Delta_f (\phi G)
\ge & ~2 \left[ \frac{(\beta-1) \nabla h}{\beta h} + \frac{\nabla \phi}{\phi} \right] \nabla(\phi G) 
- 2 G \left[ \frac{(\beta-1) \nabla h}{\beta h} + \frac{\nabla \phi}{\phi} \right] \nabla \phi \nonumber \\
&+ [2 -(2+\alpha^2)\beta-(2\beta-1)\alpha ]/(\beta  h^{\alpha +2}) \phi G^2 \nonumber\\
&+ G (\Delta_f \phi - 2(n-1)k \phi) + \mathsf{Z}_\Sigma \phi.
\end{align}

Let $q$ be a maximum point for $\phi G$ in $\mathcal{B}_{R}$. 
For the sake of establishing the estimate at $x$ in $\mathcal{B}_{R/2}$ we confine to the case $d(q) \ge 1$ noting $(\phi G)(q) >0$ . 
Now at the maximum point $q$ we have the inequalities $\Delta_f(\phi G) \le 0$ and $\nabla(\phi G) =0$. Therefore applying these to 
\eqref{Pf thm18 eq5} and rearranging the inequality we have
\begin{align} \label{eq4.2}
\left[\frac{2 -(2+\alpha^2)\beta-(2\beta-1)\alpha}{\beta} \right] \phi G^2
\le & ~2 h^{\alpha +2} G \left[ \frac{(\beta-1) \nabla h}{\beta h} + \frac{\nabla \phi}{\phi} \right] \nabla \phi \nonumber\\
&- h^{\alpha+2} G (\Delta_f \phi -2(n-1)k \phi) - h^{\alpha +2} \mathsf{Z}_\Sigma \phi.  
\end{align}

We now proceed onto bounding from above each of the terms on the right-hand side of \eqref{eq4.2}. Again, the argument proceeds by considering 
two case $d(q) \le 1$ and $d(q) \ge 1$, and so as noted above, in view of certain similarities with the proof of Theorem \ref{thm1}, we shall 
remain brief, focusing on case two only and mainly so on the differences. Towards this end, the first two terms on the right-hand side of \eqref{eq4.2} 
are seen to be bounded, directly in modulus, upon using the Cauchy-Schwarz and Young inequalities, by the expressions:  
\begin{align} \label{eq-5.16}
2 \frac{\beta-1}{\beta} h^{\alpha+1} G \langle\nabla h, \nabla \phi\rangle &\le 2 \left| \frac{\beta-1}{\beta} h^{\alpha+1} G \langle \nabla h, \nabla \phi \rangle \right| 
\le \frac{2(1-\beta)}{\beta} \phi^{3/4} G^{3/2} h^{\alpha/2+1} \frac{|\nabla \phi|}{\phi^{3/4}} \nonumber \\
& \le \frac{1-\beta}{4\beta} G^2 \phi + C(\beta) \frac{|\nabla \phi|^4}{\phi^3} h^{2 \alpha+4} 
\le \frac{1-\beta}{4\beta} G^2 \phi + \frac{C(\beta)}{R^4} h^{2\alpha+4},   
\end{align}
where we have used $\sqrt G = h^{\alpha/2} |\nabla h|$, and in much the same way, 
\begin{align} \label{eq-5.17}
2 \frac{|\nabla \phi|^2}{\phi}h^{\alpha+2} G & \le 2 \phi^{1/2} G \frac{|\nabla \phi|^2}{\phi^{3/2}} h^{\alpha +2} 
\le \frac{1-\beta}{4\beta} G^2 \phi + C(\beta) \frac{|\nabla \phi|^4}{\phi^3} h^{2\alpha +4} \nonumber\\
&\le \frac{1-\beta}{4\beta} G^2 \phi + \frac{C(\beta)}{R^4} h^{2\alpha+4}. 
\end{align}
Regarding the fourth term, i.e., the one involving $\mathsf{Z}_\Sigma$, we have upon substitution from \eqref{Z-Sigma}
\begin{align} \label{eq-5.18}
- h^{\alpha+2} \mathsf{Z}_\Sigma \phi &= - \frac{[2-(2+\alpha)\beta] \Sigma - 2u \Sigma_u}{u} h^{\alpha+2} \phi G
+ \frac{2\beta h^{2\alpha +3} \langle \nabla h, \Sigma_x \rangle \phi}{u} \nonumber \\
&\le \left[ \frac{2u \Sigma_u - [2-(2+\alpha)\beta] \Sigma}{u} \right]_+ h^{\alpha +2} \phi G 
+ 2\beta h^{3\alpha/2+3} \sqrt G \frac{|\Sigma_x|}{u} \phi \nonumber \\
&\le \frac{1-\beta}{4 \beta} G^2 \phi + C(\beta) \left[ \mathsf{T}_\Sigma^2 (u) + \mathsf{S}_\Sigma^{4/3} (u) \right] h^{2\alpha+4}, 
\end{align}
where in the last line we have written $\mathsf{T}_\Sigma (u) = \{ [2u \Sigma_u-(2-(2+\alpha)\beta) \Sigma]/u\}_+$ 
and $\mathsf{S}_\Sigma (u) = |\Sigma_x|/u$. 

Lastly, for the term involving $\Delta_f \phi$, we proceed similar to the proof of Theorem \ref{thm1}, 
where by recalling the weighted Laplacian comparison theorem as in Section \ref{sec3}, we have  
\begin{align*}
- \Delta_f \phi &  \leq C \left( \frac{1}{R^2} + \frac{[\gamma_{\Delta_f}]_+}{R} + k \right) \sqrt \phi. 
\end{align*}
Hence for the remaining term on the right-hand side of \eqref{eq4.2} upon recalling $0 \le \phi \le 1$ 
and adjusting the constant $C>0$ if necessary we can write  
\begin{align} \label{Lap comp 2 sum}
-h^{\alpha+2}G \left( \Delta_f -2(n-1)k \right) \phi 
& \leq C \left( \frac{1}{R^2} + \frac{[\gamma_{\Delta_f}]_+}{R} + k \right) h^{\alpha+2} \sqrt \phi G\nonumber \\
& \leq \frac{1-\beta}{4\beta} \phi G^2
+ C(\beta) \left( \frac{1}{R^4} + \frac{[\gamma_{\Delta_f}]_+^2}{R^2} + k^2 \right) h^{2\alpha+4}.  
\end{align}

Having now estimated each of the individual terms on the right-hand side of \eqref{eq4.2} we proceed next by substituting these back 
into the inequality and finalising the estimate. Indeed, substituting \eqref{eq-5.16}-\eqref{Lap comp 2 sum} in \eqref{eq4.2}, adding the 
expressions on the right-hand side of the former four inequalities and making note of the basic relation 
\begin{align}
\frac{2 -(2+\alpha^2)\beta-(2\beta-1)\alpha}{\beta} 
- \frac{4(1-\beta)}{4 \beta} 
&= \frac{1-(1+\alpha^2)\beta-(2\beta-1)\alpha}{\beta} \nonumber \\
&= \frac{(1+ \alpha)[1-\beta(1+\alpha)]}{\beta}, 
\end{align}
it follows that at the point $q$ we have   
\begin{align} \label{ing eq1}
\bigg[ \frac{1-(1+\alpha^2)\beta-(2\beta-1)\alpha}{\beta} \bigg] & \phi G^2 
= \bigg[ \frac{(1+ \alpha)[1-\beta(1+\alpha)]}{\beta} \bigg] \phi G^2 \\
\le&~C(\beta) \left\{ \frac{1}{R^4} + \frac{[\gamma_{\Delta_f}]_+^2}{R^2} + k^2 
+ \left[ \mathsf{T}_\Sigma^2 (u) + \mathsf{S}_\Sigma^{4/3} (u) \right] \right\} h^{2\alpha+4}. \nonumber
\end{align}

Next, by the maximality of the localised function $\phi G$ on $\mathcal{B}_R$ at $q$, we have the chain of inequalities 
[recall that $\phi \equiv 1$ on $\mathcal{B}_{R/2}$ and $0 \le \phi \le 1$ on $\mathcal{B}_R$]
\begin{equation}
\sup_{{\mathcal B}_{R/2}} G^2 = \sup_{{\mathcal B}_{R/2}} (\phi G)^2 
\le \sup_{{\mathcal B}_R} (\phi G)^2 = (\phi G)^2 (q) \le (\phi G^2)(q). 
\end{equation}

Hence combining the latter with \eqref{ing eq1} and making note of the relations 
$h=u^\beta$ and $G = h^{\alpha} |\nabla h|^2 = \beta^2 u^{(\alpha+2)\beta -2} |\nabla u|^2$ 
and introducing $s(\alpha, \beta) = \sqrt{(1+ \alpha)[1-\beta(1+\alpha)] / \beta}$ it follows that 
\begin{align*}
s(\alpha, \beta)  \frac{h^{\alpha} |\nabla h|^2}{\beta^2} 
=&~s(\alpha, \beta)
\frac{|\nabla u|^2}{u^{2-(\alpha+2)\beta}} \nonumber\\
&\le C(\beta) \left\{ \frac{1}{R^2} 
+ \frac{[\gamma_{\Delta_f}]_+}{R} + k 
+ \sup_{\mathcal{B}_{R}} \left[ \mathsf{T}_\Sigma (u) + \mathsf{S}_\Sigma^{2/3} (u) \right] \right\} 
\Big( \sup_{\mathcal{B}_{R}}  u \Big)^{(\alpha+2)\beta}. 
\end{align*}
(Note that the condition $0<\beta<1/(1+\alpha)$ guarantees that $(1+ \alpha)[1-\beta(1+\alpha)]>0$ and so in turn $s(\alpha, \beta)>0$.) 
This upon rearranging terms and taking square roots gives the desired estimate for every $x \in \mathcal{B}_{R/2}$ and so completes 
the proof. \hfill $\square$

\begin{remark} \label{Ricmf-proof-2-remark} {
A similar argument to the one outlined in Remark \ref{Ricmf-proof-1-remark} implies that subject 
to the stronger curvature assumption there, in the local estimate above the 
functional term \eqref{extra-term-Ricmf} can be removed from the right-hand 
side in \eqref{eq1.12}. See also Remark \ref{other-curvature-assuptions-remark}.
}
\end{remark}

\section{Proof of the Liouville-type result in Theorem \ref{thm Liouville}} \label{Liouville-Section}

Since here $\Sigma$ depends only on $u$ (and not $x$), by referring to \eqref{eq2.10} it is easily seen that  
$\mathsf{S}_\Sigma (u) = |\Sigma_x(x,u)|/u \equiv 0$. Moreover referring to \eqref{eq2.10-alt} and invoking 
the condition $[1-\beta(\alpha/2+1)]\Sigma(u)-u\Sigma_u(u) \ge 0$ as stated in the theorem and noting 
$u>0$ it is also seen that
\begin{equation} 
\mathsf{T}_\Sigma (u) = (2/u) [u \Sigma_u(u)-[1-\beta(1+\alpha/2)]\Sigma(u)]_+ \equiv 0.
\end{equation}
Now as $u$ is positive and bounded and we have $\beta(1+\alpha/2)>0$ it is evident that  
\begin{equation}
\sup_{M} u^{\beta(1+\alpha/2)} <\infty.
\end{equation} 
Thus in view of the assumption ${\mathscr Ric}_f(g) \ge 0$ on $M$ it follows from the global bound in Theorem 
\ref{thm18-global} with $k=0$ and after incorporating the above that  
\begin{align}
\frac{|\nabla u|}{u^{1-\beta(\alpha/2+1)}} 
&\le C \left\{ \sqrt k + \sup_{M} \left[ \mathsf{T}_\Sigma^{1/2} (u) + \mathsf{S}_\Sigma^{1/3} (u) \right] \right\}
\Big( \sup_{M} u^{\beta(\alpha/2+1)} \Big) \nonumber \\ 
&\le C \sup_{M} \left[ \mathsf{T}_\Sigma^{1/2} (u) + \mathsf{S}_\Sigma^{1/3} (u) \right] 
\Big( \sup_{M} u^{\beta(\alpha/2+1)} \Big) =0. 
\end{align}
It therefore follows that $|\nabla u| \equiv 0$ on $M$ and so $u$ is a constant. The final assertion now follows by 
substituting this constant solution $u$ back in the equation $\Delta_f u + \Sigma(u) =0$. \hfill $\square$

\section{Proof of the global Hamilton bound in Theorem \ref{global-Hamilton-Thm}} 
\label{Hamilton-Section}

\begin{lemma} \label{global-Hamilton-lemma-7}
Let $(M, g, e^{-f} dv_g)$ be a complete smooth metric measure space verifying the global curvature lower bound 
${\mathscr Ric}_f (g) \ge - \mathsf{k} g$  on $M$ with $\mathsf{k} \ge 0$. Let u be a bounded positive solution to 
\eqref{eq11} with $0<u\leq A$ and let 
\begin{align} \label{P-equation}
\mathscr{P}_\gamma[u](x)  = \gamma |\nabla u|^2/u -u \log(A/u),
\end{align}
where $\gamma$ is a non-negative constant. Then we have the elliptic differential inequality 
\begin{align}
\Delta_f \mathscr{P}_\gamma [u] 
\ge&~(1-2\gamma \mathsf{k}) \frac{|\nabla u|^2}{u} - \frac{2\gamma}{u} \langle \nabla u, \Sigma_x (x,u)\rangle \nonumber \\
& +\gamma \frac{|\nabla u|^2}{u} \left[\frac{\Sigma(x,u)}{u} - 2 \Sigma_u(x,u) \right] 
+ [\log (A/u)-1] \Sigma (x,u) 
\end{align}
\end {lemma}

\begin{proof}
Starting with the identity $\nabla (|\nabla u|^2/u) = (\nabla |\nabla u|^2)/u - (|\nabla u|^2 \nabla u)/u^2$, a further 
differentiation, upon recalling the identity $\Delta_f v = \Delta v - \langle \nabla f, \nabla v \rangle$ gives,
\begin{align}\label{eq4.43}
\Delta_f \left[ \frac{|\nabla u|^2}{u} \right] 
&= \frac{1}{u} \Delta_f |\nabla u|^2 - \frac{2}{u^2} \langle \nabla |\nabla u|^2 , \nabla u \rangle 
- \frac{|\nabla u|^2}{u^2} \Delta_f u + 2\frac{|\nabla u|^4}{u^3}.
\end{align}

Now referring to the right-hand side and by using the weighted Bochner-Weitzenb\"ock formula 
$\Delta_f |\nabla u|^2 = 2 |\nabla ^2 u|^2 + 2\langle \nabla u , \nabla \Delta_f u \rangle + 2 {\mathscr Ric}_f (\nabla u, \nabla u)$ we can write 
\begin{align}\label{eq3.45}
\Delta_f \left[ \frac{|\nabla u|^2}{u} \right] =&~ \frac{2}{u} [ |\nabla ^2 u|^2 
+ \langle \nabla u , \nabla \Delta_f u \rangle +  {\mathscr Ric}_f (\nabla u, \nabla u)]\nonumber\\
& - \frac{2}{u^2} \langle \nabla |\nabla u|^2 , \nabla u \rangle - \frac{|\nabla u|^2}{u^2} \Delta_f u + 2\frac{|\nabla u|^4}{u^3} \\ 
=&~ \frac{2}{u} \left| \nabla ^2 u - \frac{\nabla u \otimes \nabla u }{u} \right| ^2
+\frac{2}{u} {\mathscr Ric}_f (\nabla u, \nabla u) \nonumber \\
&-\frac{2\langle \nabla u,\nabla \Sigma(x,u) \rangle }{u} + \frac{|\nabla u|^2}{u^2} \Sigma(x,u), \nonumber 
\end{align}
where we have made use of \eqref{eq11} and 
\begin{align}
|\nabla ^2 u|^2 - \frac{1}{u} \langle \nabla |\nabla u|^2 , \nabla u \rangle 
+ \frac{|\nabla u|^4}{u^2} = \left| \nabla ^2 u - \frac{\nabla u \otimes \nabla u}{u} \right| ^2.
\end{align}

Likewise a straightforward differentiation gives $\nabla [u \log (A/u)] = [\log (A/u) -1] \nabla u$ and so calculating the $f$-Laplacian we have 
\begin{align}\label{eq4.50}
\Delta_f [u \log (A/u)] = [\log (A/u) -1] \Delta_f u - |\nabla u|^2/u.  
\end{align}
As a result by putting the two segments together and noting $\Delta_f u = -\Sigma(x,u)$ it is seen that 
\begin{align}\label{eq4.51}
\Delta_f \mathscr{P}_\gamma[u] =&~\Delta_f [\gamma |\nabla u|^2/u - u \log(A/u)] \nonumber\\
=& \frac{2\gamma}{u} \left| \nabla ^2 u - \frac{\nabla u \otimes \nabla u}{u} \right| ^2
+\frac{2\gamma}{u} {\mathscr Ric}_f (\nabla u, \nabla u) 
+ \gamma\frac{|\nabla u|^2}{u^2}\Sigma(x,u) \nonumber \\ 
& - \frac{2\gamma}{u} \langle \nabla u,\nabla \Sigma (x,u) \rangle 
+ |\nabla u|^2/u + [\log (A/u) -1] \Sigma(x,u).
\end{align}
 Next by virtue of the inequalities 
\begin{align}
\left| \nabla ^2 u - \frac{\nabla u \otimes \nabla u }{u} \right| ^2 \ge 0, \qquad \text{and} \qquad  {\mathscr Ric}_f (g) \ge - \mathsf{k} g, 
\end{align}
it follows from \eqref{eq4.51} that
\begin{align}
\Delta_f \mathscr{P}_\gamma[u] \ge&~(1-2\gamma \mathsf{k}) \frac{|\nabla u|^2}{u} 
-\frac{2\gamma}{u} \langle \nabla u,\nabla \Sigma (x,u) \rangle \nonumber\\
&+ \left[ \gamma \frac{|\nabla u|^2}{u^2} + [\log (A/u) -1] \right] \Sigma(x,u) \\
\ge&~(1-2\gamma \mathsf{k}) \frac{|\nabla u|^2}{u} 
- \frac{2\gamma}{u} \langle \nabla u, \Sigma_x (x,u)\rangle \nonumber \\
& +\gamma \frac{|\nabla u|^2}{u} \left[\frac{\Sigma(x,u)}{u} - 2 \Sigma_u(x,u) \right] 
+ [\log (A/u)-1] \Sigma (x,u),  
\end{align}
which is the desired conclusion. 
\end{proof}

 \qquad \\
{\it Proof of Theorem $\ref{global-Hamilton-Thm}$.} 
If $u>0$ is constant then \eqref{eq7.11} is trivially true so we assume $u$ is non-constant. Now since for any 
$\sigma \ge e$ we have $u \le \sigma A$, an application of Lemma \ref{global-Hamilton-lemma-7} above, 
upon replacing $A$ with $\sigma A$ and noting $\Sigma_x \equiv 0$ gives, 
\begin{align}
\Delta_f \mathscr{P}_\gamma [u] 
\ge&~(1-2\gamma \mathsf{k}) \frac{|\nabla u|^2}{u} + [\log (\sigma A/u)-1] \Sigma (u) \nonumber \\
& +\gamma \frac{|\nabla u|^2}{u} \left[\frac{\Sigma(u)}{u} - 2 \Sigma'(u) \right]  \nonumber \\
\ge&~(1-2\gamma \mathsf{k}) \frac{|\nabla u|^2}{u} + [\log (\sigma/e) + \log (A/u)] \Sigma (u)  \nonumber \\
& +\gamma \frac{|\nabla u|^2}{u} \left[\frac{\Sigma(u)}{u} - 2 \Sigma'(u) \right] \ge 0,
\end{align}
subject to setting $\gamma=1/(2 \mathsf{k})$ where we have used $\Sigma(u) \ge 0$ and 
$\Sigma(u) - 2 \Sigma'(u)/u \ge 0$. Since $M$ is closed an application of the maximum principle now gives 
$\mathscr{P}_\gamma[u] \le 0$ or else that $\mathscr{P}_\gamma[u]$ is constant. However since $u$ is 
non-constant and $\sigma \ge e$ is arbitrary, the function 
$\mathscr{P}_\gamma[u] = \gamma |\nabla u|^2/u - u [\log (A /u) + \log \sigma]$ 
is constant for at most one $\sigma \ge e$. Thus we can write 
\begin{equation} 
|\nabla \log u|^2 = |\nabla u|^2/u^2 \le 2 \mathsf{k} [\log (A /u) + \log \sigma]
\end{equation}
and so by passing to the limit $\sigma \searrow e$ and noting that $M$ is compact we obtain \eqref{eq7.11}. 
Next to prove \eqref{eq16} set $\mathcal{Z}(x) = \log [(Ae)/u]=1+\log(A/u)$. Then a straightforward calculation 
and making use of \eqref{eq7.11} gives
\begin{equation}
\left| \nabla \sqrt{\mathcal{Z}} \right| = \left| \frac{\nabla u/u}{2 \sqrt{\mathcal{Z}}} \right| 
\le \frac{\sqrt{2 \mathsf{k} {\mathcal Z}}}{2 \sqrt{\mathcal Z}} = \sqrt{\mathsf{k}/2}.
\end{equation}
Integrating the above along a minimising geodesic joining a pair of points $p, q$ in $M$ then gives 
\begin{equation}
\sqrt{\log (Ae/u(q))} \le \sqrt{\log (Ae/u(p))} + d(p,q) \sqrt{{\mathsf{k}/2}}.
\end{equation} 
For any $\varepsilon >0$ thus 
$\log (Ae/u(q)) \le (1+\varepsilon) [ \log (Ae/u(p)) + d^2(p,q) {\mathsf k}/(2\varepsilon)]$. 
Exponentiating and rearranging yields the desired inequality \eqref{eq16}. 
\hfill $\square$


\qquad \\
{\bf Acknowledgement.} The authors gratefully acknowledge support from EPSRC. They also wish to thank the anonymous 
referee for a careful reading of the manuscript and useful comments.


\begin{thebibliography}{99}

\bibitem{AM} E.~Acerbi, R.~Mingione, {\it Gradient estimates for a class of parabolic systems}, Duke Math. J., 136, (2007), 285-320.

\bibitem{BD}{D. Bakry}, {\it L'hypercontractivit\'e et son utilisation en th\'eorie des semigroupes}, in: Lecture 
Notes in Math., {\bf 1581}, Springer-Verlag, Berlin/New York, (1994), 1--114. 
					
\bibitem{BE}{D. Bakry, M. \'Emery}, {\it Diffusions hypercontractives} In: Az\'ma J., Yor M. (eds) 
S\'eminaire de Probabilit\'es XIX 1983/84. Lecture Notes in Mathematics, {\bf 1123}, Springer, Berlin, Heidelberg.

\bibitem{Bak}{D. Bakry, I. Gentil, M. Ledoux}, {\it Analysis and Geometry of Markov Diffusion Operators}, 
A Series of Comprehensive Studies in Mathematics, {\bf 348}, Springer, 2012.

\bibitem{Bid}{M.F. Biduat-V\`eron, L. Ver\`on}, {\it Nonlinear elliptic equations on compact Riemannian manifolds 
and asymptotics of the Emden equations}, Invent. Math., 106, (1991), 489--539.



\bibitem{Bri}{K. Brighton}, {\it A Liouville-type theorem for smooth metric measure spaces}, 
J. Geom. Anal., 23 (2013), 562--570.

\bibitem{CGS} L.A.~Cafarreli, B.~Gidas, J.~Spruck, {\it Asymptotic symmetry and local behavior of semilinear 
elliptic equations with critical Sobolev growth}, Comm. Pure Appl. Math., 42, (1989), 271--297. 

\bibitem{Calabi} E.~Calabi, {\it An extension of Hopf maximum principle with application to Riemannian geometry}, 
Duke Math. J., 25, (1958), 45--56.  

\bibitem{Cao} H.~Cao, {\it Recent progress on Ricci solitons}, In: Recent Advances in Geometric Analysis, 
Advanced Lectures in Mathematics (ALM), Vol.~11, 1--38. International Press, (2010)

\bibitem{Case} J.~Case, {\it A Yamabe problem on smooth metric measure spaces}, J. Diff. Geom., 
101 (2015), 467--505.

\bibitem{ChYau} S.Y.~Cheng, S.T.~Yau, {\it Differential equations on Riemannian manifolds and their geometric applications}, 
Comm. Pure Appl. Math., 28, (1975), 333--354. 

\bibitem{CBY}{Y. Choquet-Bruhat}, {\it General Relativity and the Einstein Equations}, Oxford Mathematical Monographs, 
OUP, 2009.

\bibitem{Chow}{B. Chow, P. Lu, L. Nei}, {\it Hamilton's Ricci Flow}, Graduate Studies in Mathematics {\bf 77}, AMS, 2006.	

\bibitem{Dung}{N.T.~Dung, N.N.~Khanh, Q.A.~Ng\^o}, {\it Gradient estimates for $f$-heat equations driven by 
Lichnerowicz's equation on complete smooth metric measure spaces}, Manuscripta Math., 155, (2018), 471--501.

\bibitem{GKS} M.~Ghergu, S.~Kim, H.~Shahgholian, {\it Exact behaviour around isolated singularity for semilinear 
elliptic equations with a log-type nonlinearity}, Adv. Nonlinear Anal., 8 (2019), 995--1003.
 
\bibitem{Giaq} M.~Giaquinta, {\it Multiple Integrals in the Calculus of Variations and Nonlinear Elliptic Systems}, 
Annals of Mathematics Studies, Vol.~{\bf 105}, Princeton University Press, 1983. 
 
\bibitem{GidSp} B.~Gidas, J.~Spruck, {\it Global and local behaviour of positive solutioins of nonlinear elliptic equations}, 
Comm. Pure Appl. Math., 34 (1981), 525--598. 

\bibitem{Gr}{A. Grigor'yan}, {\it Heat kernel analysis on manifolds}, Studies in Advanced Mathematics, AMS, 2013. 

\bibitem{Gross} L.~Gross, {\it Logarithmic Sobolev inequalities}, Amer. J. Math., 97, (1976), 1061--1083. 

  
\bibitem{Ha93}{R. Hamilton},  {\it A matrix Harnack estimate for heat equation}, Comm. Anal. Geom., (1993), 113--126.	 
 
\bibitem{Ham} R.~Hamilton, {\it The formation of singularities in the Ricci flow}, Surv. Diff. Geom., 2,  (1995), 7--136. 
 
\bibitem{KT}{J.~Kristensen, A.~Taheri}, {\it Partial regularity of strong local minimizers in the calculus of variations}, 
Arch. Rational Mech. Anal., 170, (2003), 63--89.
  
\bibitem{Lee} J.M.~Lee, T.H.~Parker, {\it The Yamabe problem}, Bull. Amer. Math. Soc., 17, (1987), 37--91.
 
\bibitem{LiJ91} {J.~Li}, {\it Gradient estimates and Harnack inequalities for nonlinear parabolic and nonlinear elliptic 
equations on Riemannian manifolds}, J. Funct. Anal., 100, (1991), 233--256. 

\bibitem{LY86}{P. Li, S.T. Yau}, {\it On the parabolic kernel of Schr\"odinger operator}, Acta Math., 156 (1986), 153--201. 
 
\bibitem{LiPTam} P.~Li, L.F.~Tam, D.G.~Yang, {\it On the elliptic equation $\Delta u + k u - K u^p=0$ on complete 
Riemannian manifolds and their geometric applications}, Trans. Amer. Math. Soc., 350, (1998), 1045--1078.     
 
\bibitem{LiP-book}{P. Li}, {\it Geometric Analysis}, Cambridge Studies in Advanced Mathematics, {\bf 134}, CUP, 2012.

\bibitem{LiX}{X.D. Li}, {\it Liouville theorems for symmetric diffusion operators on complete Riemannian manifolds}, 
J. Math. Pures Appl., 84 (2005) 1295--1361.

\bibitem{Lott}{J. Lott}, {\it Some geometric properties of the Bakry-\'Emery Ricci tensor}, Comment. Math. Helv., 78 (2003), 865--883. 

\bibitem{Ma} L.~Ma, {\it Gradient estimates for a simple elliptic equation on complete noncompact Riemannian manifolds}, 
J. Funct. Anal., 241 (2006), 374--382. 

\bibitem{Mast} P.~Mastrolia, M.~Rigoli, A.G.~Setti, {\it Yamabe Type Equations on Complete Non-compact Manifolds}, 
Springer, Basel 2012.


\bibitem{PG}{G. Perelman}, {\it The entropy formula for the Ricci Flow and its geometric application}, arXiv: math. DG/0211159v1 (2002).
 
\bibitem{Qian1} Z.~Qian, {\it On conservation of probability and the Feller property}, Ann. Prob., 24 (1996), 280--292.  

 
\bibitem{Ruan}{Q.H.~Ruan} {\it Elliptic type gradient estimates for Schr\"odinger equations on noncompact manifolds}, 
Bull. Lond. Math. Soc., 39 (2007), 982--988. 
 
\bibitem{SchYau} R.~Schoen, S.T.~Yau, {\it Lectures on Differential Geometry}, International Press, 1994. 
 
 
\bibitem{SZ} {P. Souplet, Q.S. Zhang} 
{\it Sharp gradient estimate and Yau's Liouville theorem for the heat equation on noncompact manifolds}, 
Bull. Lond. Math. Soc., 38 (2006), 1045--1053. 
 
\bibitem{Taheri-book-one} A.~Taheri, 
{\it Function Spaces and Partial Differential Equations}, Vol.~{\bf I}, Oxford Lecture Series in Mathematics and its Applications, {\bf 40}, OUP, 2015.   

\bibitem{Taheri-book-two} A. Taheri, 
{\it Function Spaces and Partial Differential Equations}, Vol.~{\bf II}, Oxford Lecture Series in Mathematics and its Applications, {\bf 41}, OUP, 2015. 

\bibitem{Taheri-GE-1} A.~Taheri, {\it Liouville theorems and elliptic gradient estimates for a nonlinear parabolic equation involving the Witten Laplacian}, 
Published online in: Adv. Calc. Var., De Gruyter, 2021.

\bibitem{Taheri-GE-2} A.~Taheri, {\it Gradient estimates for a weighted $\Gamma$-nonlinear parabolic equation coupled with a super Perelman-Ricci flow 
and implications}, Published online in: Potential Anal., Springer, 2021. 

\bibitem{TVah} A.~Taheri, V.~Vahidifar, {\it On multiple solutions to a family of nonlinear elliptic systems in divergence form combined 
with an incompressibility constraint}, Nonlinear Anal., 221, 2022. 

\bibitem{TVahGrad} A.~Taheri, V.~Vahidifar, {\it Gradient estimates for nonlinear elliptic equations involving the Witten Laplacian on smooth 
metric measure spaces and implications}, Adv. Nonlinear Anal., 12, De Gruyter, 2023.

\bibitem{TahVah} A.~Taheri, V.~Vahidifar, {\it Gradient estimates for a nonlinear parabolic equation on smooth metric measure spaces 
with evolving metrics and potentials}, Nonlinear Anal., 232, 2023.

\bibitem{VC} {C. Villani}, {\it Optimal transport: Old and New}, A Series of Comprehensive Studies in Mathematics, {\bf 338}, Springer, 2008.
 
 
\bibitem{Wang}{F.Y.~Wang}, {\it Analysis for Diffusion Processes on Riemannian Manifolds}, Advanced Series on Statistical Science $\&$ Probability, 
Vol.~{\bf 18}, World Scientific, 2013.
 
\bibitem{[WeW09]}{G.~Wei, W.~Wylie}, {\it Comparison geometry for the Bakry-\'Emery Ricci tensor}, J. Diff. Geom., 83 (2009), 377--405.


\bibitem{Wu15}{J.Y.~Wu}, {\it Elliptic gradient estimates for a weighted heat equation and applications}, Math. Z., 280 (2015), 451--468.

\bibitem{Wu18}{J.Y.~Wu}, {\it Gradient estimates for a nonlinear parabolic equation and Liouville theorems}, Manuscript Math., 159 (2019), 511--547. 


\bibitem{WuLM} L.M.~Wu, {\it Uniqueness of Nelson's diffusion}, Probability Theory and Related Fields, 114, (1999), 549--585.

\bibitem{YYY} Y.Y.~Yang, {\it Gradient estimates for a nonlinear parabolic equation on Riemannian manifolds}, Proc. Amer. Math. Soc., 136 (2008), 4095--4102.

\bibitem{Yang} Y.Y.~Yang, {\it Gradient estimates for the equation $\Delta u + c u^{-\alpha}=0$ on Riemannian manifolds}, Acta Math. Sin. 26, (2010), 1177-1182.

\bibitem{YauH} S.T.~Yau, {\it Harmonic functions on complete Riemannian manifolds}, Comm. Pure Appl. Math., 28 (1975), 201-228.

\bibitem{ZhMa} J.~Zhang, B.~Ma, {\it Gradient estimates for a nonlinear equation $\Delta_f u + cu^{-\alpha} =0$ on complete noncompact manifolds}, 
Comm. Math., 19, (2011), 73-84.


\bibitem{Zhang} Q.S.~Zhang, {\it Sobolev inequalities, heat kernels under Ricci flow and the Poincar\'e conjecture}, CRC Press, 2011. 





 
 \end{thebibliography}
\end{document}